\documentclass{amsart}
\usepackage{FrankTall}
\usepackage[usenames,dvipsnames]{xcolor}
\usepackage{xspace}
\usepackage{accents}
\usepackage{dirtytalk}
\usepackage[normalem]{ulem}

\theoremstyle{plain}
\newtheorem{ithm}{Theorem}

\newcommand{\R}{\mathbb{R}}
\newcommand{\w}{\mathbb{\omega}}
\newcommand{\script}[1]{\mathcal{#1}}
\newcommand{\restr}{\upharpoonright}

\newcommand{\ignore}[1]{}
\newcommand{\VL}{\ensuremath{\mathsf{V=L}}}
\newcommand{\MA}{\ensuremath{\mathsf{MA}}}
\newcommand{\bSigma}{\mathbf{\Sigma}}
\newcommand{\bPi}{\mathbf{\Pi}}
\newcommand{\bDelta}{\mathbf{\Delta}}

\newcommand{\comp}[1]{\ensuremath{\overline{#1}}} 

\newcommand{\Mod}{\mathrm{Mod}} 
\newcommand{\opgm}{\Game} 
\newcommand{\pair}[2]{\langle #1, #2 \rangle}

\newcommand{\stp}{\,\,\,}

\begin{document}
\title[An undecidable extension of Morley's theorem]{An undecidable extension of Morley's theorem on the number of countable models}
\author[C. J. Eagle]{Christopher J. Eagle${}^1$} 
\address[C. J. Eagle]{University of Victoria, Department of Mathematics and Statistics, PO BOX 1700 STN CSC, Victoria, British Columbia, Canada, V8W 2Y2}%
\email{eaglec@uvic.ca}
\urladdr{http://www.math.uvic.ca/~eaglec}
\thanks{${}^1$ Supported by NSERC Discovery Grant RGPIN-2021-02459}

\author[C. Hamel]{Clovis Hamel${}^2$}
\address[C. Hamel]{University of Toronto, Department of Mathematics, 40 St. George St., Toronto, Ontario, Canada M5S 2E4}
\email{chamel@math.toronto.edu}
\urladdr{}
\thanks{${}^2$ Supported by an NSERC Vanier Scholarship.}

\author[S. M\"uller]{Sandra M\"uller${}^3$}
\address[S. M\"uller]{Institut f\"ur Diskrete Mathematik und Geometrie, TU Wien, Wiedner Hauptstra{\ss}e 8-10/104, 1040 Wien, Austria}
\email{sandra.mueller@tuwien.ac.at}
\urladdr{https://dmg.tuwien.ac.at/sandramueller/}
\thanks{${}^3$ Supported by L'OR\'{E}AL Austria, in collaboration with the Austrian UNESCO Commission and in cooperation with the Austrian Academy of Sciences - Fellowship \emph{Determinacy and Large Cardinals} and the Austrian Science Fund (FWF) under Elise Richter grant number V844, international project number I6087, and START grant number Y1498}

\author[F. D. Tall]{Franklin D. Tall${}^4$}
\address[F. D. Tall]{University of Toronto, Department of Mathematics, 40 St. George St., Toronto, Ontario, Canada M5S 2E4}
\email{f.tall@utoronto.ca}
\urladdr{http://www.math.toronto.edu/tall/}
\thanks{${}^4$ Supported by NSERC grant RGPIN-2016-06319}
\date{\today}



\begin{abstract}
    We show that Morley's theorem on the number of countable models of a countable first-order theory becomes an undecidable statement when extended to second-order logic.  More generally, we calculate the number of equivalence classes of equivalence relations obtained by countable intersections of projective sets in several models of set theory.  Our methods include random and Cohen forcing, Woodin cardinals and Inner Model Theory.
\end{abstract}
\maketitle
\section{Introduction}
\renewcommand\theithm{\Alph{ithm}}

\emph{Vaught's Conjecture}, which asserts that a countable first-order theory must have either at most countably many or exactly $2^{\aleph_0}$ many non-isomorphic countable models, is one of the most important problems in Model Theory.  While the question itself is model-theoretic, it is known to have deep connections to both Descriptive Set Theory and Topological Dynamics.  Since Vaught's original paper \cite{Vaught}, Vaught's Conjecture has been verified for a number of classes of theories, such as theories of trees \cite{Steel}, $\omega$-stable theories \cite{ShelahHarringtonMakkai}, o-minimal theories \cite{Mayer}, and varieties (in the sense of universal algebra) \cite{HartStarchenkoValeriote}, among others.

A strong positive result about Vaught's Conjecture that applies to all first-order theories is a result of the late Michael Morley \cite{Morley1970}, which states that the number of isomorphism classes of countable models of a countable first-order theory is always either at most $\aleph_1$ or exactly $2^{\aleph_0}$.  In this form Morley's Theorem requires no proof at all if we are in a universe of set theory where the continuum hypothesis holds.  However, not long after Morley's paper appeared it was noticed that this result can be improved.  Under a standard identification between countable models and elements of $2^\omega$ (described in Section \ref{sec:Preliminaries} below), we have the following strengthening of Morley's result:

\begin{thm*}[Absolute Morley Theorem]
Let $T$ be a first-order theory (or, more generally, a sentence of $L_{\omega_1, \omega}$) in a countable signature.  Then either $T$ has at most $\aleph_1$ isomorphism classes of countable models, or there is a perfect set of non-isomorphic countable models of $T$.
\end{thm*}

The previous theorem does not appear in Morley's paper, but the ideas are there for its proof. The theorem immediately follows from an important Descriptive Set Theory result by Burgess {\cite[Corollary 2]{Burgess2}}:

\begin{thm*}[Burgess]
Let $E$ be a $\bSigma^1_1$ equivalence relation on $\R$. If there is no perfect set of pairwise inequivalent reals, then there are at most $\aleph_1$ equivalence classes.
\end{thm*}

The intersection of Descriptive Set Theory and Model Theory has become a widely studied topic, to such extent that there are versions of Vaught's Conjecture that completely fall in the domain of Descriptive Set Theory, as they concern equivalence relations more general than isomorphism between countable models (see \cite{Gao}). The main idea behind the connection between these two disciplines is that countable models of a second-order theory can be \textit{coded} as reals, usually as elements of $2^\omega$. It is then easy to prove that the isomorphism relation (which can be formulated as the existence of a certain function) is a $\bSigma_1^1$ equivalence relation. Finally, the proof of the Absolute Morley Theorem follows from the fact proved by Morley \cite{Morley1970} that the set of countable models of a theory in a countable fragment of $L_{\omega_1, \omega}$ is Borel.  

In this paper we concern ourselves with logics that are stronger than first-order logic, especially second-order logic and fragments thereof.  Vaught's Conjecture is false for second-order logic if the continuum hypothesis fails; one easy counterexample is that one can express in second-order logic that a linear order is a well-order, and hence there is a second-order theory whose countable models are (up to isomorphism) exactly the countable ordinals.  In this paper we investigate versions of the Absolute Morley Theorem for second-order logic.  Our main results are the following:

\begin{ithm}\label{thm:B}
Force over $L$ by first adding $\aleph_2$ Cohen reals and then adding $\aleph_3$ random reals.  In the resulting universe of set theory, there is a second-order theory $T$ in a countable signature such that the number of non-isomorphic countable models of $T$ is exactly $\aleph_2$, while $2^{\aleph_0} = \aleph_3$.
\end{ithm}

\begin{ithm}\label{thm:A}
Beginning with a supercompact cardinal, carry out the standard forcing iteration for producing a model of the Proper Forcing Axiom.  In the resulting universe of set theory, if $T$ is a second-order theory in a countable signature, then either $T$ has at most $\aleph_1$ isomorphism classes of countable models, or there is a perfect set of non-isomorphic models of $T$.
\end{ithm}

These two results together show that, modulo the consistency of the existence of a supercompact cardinal, the extension of the absolute version of Morley's Theorem from first-order logic to second-order logic is undecidable.

A proof of Theorem \ref{thm:A} is implicit in results of Foreman and Magidor \cite{FM95}, and in fact produces a similar conclusion not just for isomorphism classes of countable models of a second-order theory, but in fact for any equivalence relation in $L(\mathbb{R})$.  This suggests the possibility of reducing the large cardinal strength below a supercompact cardinal for obtaining our desired result.  Moreover, the full strength of second-order logic is significantly more than is needed for expressing most theories of interest in mathematics.  We therefore also wish to consider the required large cardinal strength for fragments of second-order logic with a prescribed maximum quantifier complexity (such as existential second-order logic).  We obtain the following:

\begin{ithm}\label{thm:C}
If there are infinitely many Woodin cardinals, then there is a model of set theory in which the Absolute Morley Theorem holds for second-order theories in countable signatures.
\end{ithm}

Although we believe our use of large cardinals is necessary, we don't have a proof of this.

\begin{prob}
Prove that large cardinals are necessary to prove the consistency of the conclusion
of Theorem \ref{thm:C}.
\end{prob}

The remainder of the paper is structured as follows.  In Section \ref{sec:Preliminaries} we briefly review the syntax and semantics of second-order logic, and also the method by which countable structures can be encoded as elements of $2^{\omega}$.  We also consider the descriptive complexity of the set of models of a second-order theory.  Next, in Section \ref{sec:ConsistentFailure} we prove Theorem \ref{thm:B}, and in Section \ref{sec:Consistency} we sketch a proof of a weaker version of Theorem \ref{thm:A}.  In Section \ref{sec:Extension} we develop finer results concerning the consistency of the Absolute Morley Theorem for fragments of second-order logic.  In Section \ref{sec:OtherLogics} we investigate Absolute Morley for certain game and partially ordered quantifiers.  In Section \ref{sec:Woodin} we prove Theorems \ref{thm:A} and \ref{thm:C}.

\subsection*{Acknowledgements}
We are grateful to Professor Magidor for explaining the proof of Theorem \ref{thm:ForemanMagidor} to the fourth author. We are also very grateful to the referee for a very careful reading and for catching several errors and suggesting many improvements.

\section{Preliminaries}\label{sec:Preliminaries}
\subsection{Coding countable structures as reals}
All of the logics we will consider make use of the same notion of \emph{structures} as is used in first-order logic.  For simplicity we will consider only relational structures.  Structures with function and constant symbols can be easily incorporated into this framework by coding the functions and constants as relations.  Thus, for us, a \emph{signature} consists of a collection of relation symbols, each with a specified \emph{arity}.  If $S$ is a signature, then the \emph{$S$-structures} are defined exactly as in first-order logic.  

To use tools from Descriptive Set Theory we need to code countable structures as elements of a Polish space.  The method for doing so is standard (see for example \cite[Section 11.3]{Gao}), but we review it here for the convenience of the reader.  

Let $S = \{R_i\}_{i \in I}$ be a signature, where each $R_i$ is a relation symbol of arity $n_i$, and $I$ is countable.  Suppose that $\mathcal{M}$ is a countable $S$-structure.  Up to isomorphism, we may assume that the underlying set of $\mathcal{M}$ is $\omega$.  For each $i \in I$ the interpretation $R_i^{\mathcal{M}}$ of $R_i$ is a subset of $\omega^{n_i}$, and so we identify $R_i$ with an element of $2^{\omega^{n_i}}$.  As the structure $\mathcal{M}$ is completely determined by the interpretations of each of the relation symbols, we may identify $\mathcal{M}$ with an element of $\prod_{i \in I}2^{\omega^{n_i}}$.  This identification provides a bijective map from the collection of $S$-structures with universe $\omega$ to $\prod_{i \in I}2^{\omega^{n_i}}$.  We thus view the Cantor space $\prod_{i \in I}2^{\omega^{n_i}}$ as being the space of countable $S$-structures, and define $\op{Mod}_S = \prod_{i \in I}2^{\omega^{n_i}}$.

If $\sigma$ is a sentence of some logic (such as first-order logic, second-order logic, $L_{\omega_1, \omega}$, or, more generally a model-theoretic logic as defined in \cite[Definition 1.1.1]{Ebbinghaus}), we define $\op{Mod}_S(\sigma) = \{\mathcal{M} \in \op{Mod}_S : \mathcal{M} \models \sigma\}$ (if $S$ is clear from context we may omit it).  In the case where $\sigma$ is a sentence of $L_{\omega_1, \omega}$, the set $\op{Mod}_S(\sigma)$ is a Borel subset of $\op{Mod}_S$, and moreover every isomorphism-invariant Borel subset of $\op{Mod}_S$ is of the form $\op{Mod}_S(\sigma)$ for some $L_{\omega_1, \omega}$ sentence $\sigma$ \cite[Theorem 11.3.6]{Gao}.

We are now prepared to define the main equivalence relations we study in this paper.

\begin{defi}
Let $S$ be a countable signature, and let $T$ be an $S$-theory of some logic.  The equivalence relation of \emph{isomorphism of models of $T$} is the equivalence relation $\cong_T$ on $\op{Mod}_S$ defined by declaring that $\mathcal{M} \cong_T \mathcal{N}$ if and only if either neither of the two structures is a model of $T$, or $\mathcal{M} \cong \mathcal{N}$.   For a single sentence $\sigma$, we write $\cong_\sigma$ instead of $\cong_{\{\sigma\}}$.
\end{defi}
The equivalence classes of $\cong_T$ are thus one class for each isomorphism class of models of $T$, together with one additional class containing all elements of $\op{Mod}_S \setminus \op{Mod}_S(T)$.  We are only interested in the classes corresponding to isomorphism types of models of $T$; the following lemma (which we often use without explicit mention) allows us to move from a perfect set of $\cong_T$-inequivalent structures to a perfect set of non-isomorphic models of $T$.  Recall that a set of reals $A$ has the \emph{perfect set property} if either $A$ is countable or $A$ includes a non-empty perfect set.

\begin{lem}
Let $X$ be a Polish space, and let $A \subseteq X$ be a non-empty perfect set.  For every $x \in X$, $A \setminus \{x\}$ includes a non-empty perfect set.
\end{lem}
\begin{proof}
Perfect sets are closed; in a Polish space they are therefore $G_\delta$.  Therefore the set $A \setminus \{x\}$ is a Borel set in $X$.  Since Borel sets have the perfect set property and $A \setminus \{x\}$ is uncountable, $A \setminus \{x\}$ includes a non-empty perfect set.
\end{proof}

The following proposition follows immediately from the definition of $\cong_T$, and will be key for us later.

\begin{prop}
Let $S$ be a countable signature, and let $T$ be an $S$-theory of some logic.  The descriptive set-theoretic complexity of the equivalence relation $\cong_T$ is at most the minimum projective pointclass that includes both $\bSigma_1^1$ and the complexity of the complement of $\op{Mod}_S(T)$.
\end{prop}

We end this section with a standard definition that will be used many times in what follows.

\begin{defi}
An equivalence relation $E$ on a Polish space $X$ is \textit{thin} if there is no perfect set of pairwise $E$-inequivalent elements of $X$.
\end{defi}

In particular, many of the results we are interested in are about counting the number of equivalence classes of $\cong_T$ in cases where that relation is thin.

\subsection{Second-order logic}
The logics we study in this paper are all closely related to second-order logic, so we include a review of that logic here.  The reader familiar with the ``full" semantics of second-order logic can safely skip this section, while the reader interested in learning more about second-order logic should consult \cite{Vaananen} or \cite{Vaeaenaenen2023}.

On the syntactic side the definitions closely follow the corresponding recursive definitions for first-order logic, but with additional clauses describing the use of second-order variables. Indeed, the first significant difference from first-order logic is that our formulas will allow several kinds of variables, specifically:
\begin{itemize}
    \item Variables to represent individual elements of structures.  (These are the \emph{first-order variables}).
    \item For each $n$, variables to represent sets of $n$-tuples of elements.  (These are the \emph{$n$-ary relation variables}).
\end{itemize}
The \emph{atomic second-order $S$-formulas} are defined recursively:
\begin{itemize}
    \item If $x$ and $x'$ are first-order variables, then $x=x'$ is an atomic $S$-formula.
    \item If $x_1, \ldots, x_n$ are first-order variables, and $R$ is an $n$-ary relation symbol in $S$, then $R(x_1, \ldots, x_n)$ is an atomic $S$-formula.
    \item If $x_1, \ldots, x_n$ are first-order variables, and $U$ is an $n$-ary relation variable, then $U(x_1, \ldots, x_n)$ is an atomic $S$-formula.
\end{itemize}
The \emph{second-order $S$-formulas} are defined recursively:
\begin{itemize}
    \item Atomic second-order $S$-formulas are second-order $S$-formulas.
    \item The second-order $S$-formulas are closed under conjunction, disjunction, and negation.
    \item If $\phi$ is an $S$-formula and $x$ is a first-order variable then $(\exists x) \phi$ and $(\forall x) \phi$ are second-order $S$-formulas.
    \item If $\phi$ is an $S$-formula and $U$ is a relation variable then $(\exists U) \phi$ and $(\forall U) \phi$ are second-order $S$-formulas.
\end{itemize}

In cases where it helps with clarity, we sometimes indicate that a quantifier is a second-order quantifier by adding a superscript $1$; thus the notation $\exists^1$ is sometimes used for emphasis when we are using existential quantification over a second-order variable. Analogously, $\forall^1$ is used to denote a second-order universal quantifier. 

Finally, we come to the second-order satisfaction relation.  We use this relation with the \emph{full semantics}.  These semantics are defined as for first-order logic, with the following additions.  Suppose that $\phi(P)$ is a second-order $S$-formula, where $P$ is a relation variable ($\phi$ may have other variables that are not displayed).
\begin{itemize}
    \item If $P$ is an $n$-ary relation variable and $A \subseteq M^n$, then $\mathcal{M} \models \phi(A)$ if and only if $(\mathcal{M}, A) \models \phi$, where $(\mathcal{M}, A)$ is the expanded structure obtained by interpreting $P$ as $A$.
    \item $\mathcal{M} \models (\exists P) \phi(P)$ if and only if there is some $A \subseteq M^n$ such that $\mathcal{M} \models \phi(A)$.  The definition for the second-order universal quantifier is similar.
\end{itemize}

While we will not be directly dealing with the meta-mathematical properties of second-order logic in this paper, we emphasize to the reader that second-order model theory is significantly different from first-order model theory.  In particular, both the compactness and L\"owenheim-Skolem theorems fail for second-order logic.  

In Section \ref{sec:Extension} we will need to consider restricted classes of second-order formulas.  For that purpose the following fact (see \cite[Section 4]{Vaananen}) is very useful:
\begin{fact}
Every second-order formula is equivalent to a second-order formula in \emph{prenex normal form}; that is, to a formula where all quantifiers appear at the beginning of the formula, and all second-order quantifiers precede all first-order quantifiers.
\end{fact}
In light of this fact, when discussing second-order formulas in general we will assume they are already written in prenex normal form.  To avoid overlap of terminology, we will say that a second-order formula is $\forall_n$ if it is equivalent to a prenex formula that begins with a second-order universal quantifier and has a total of $n$ blocks of quantifiers, and likewise a formula is $\exists_n$ if it is equivalent to a prenex formula that begins with a second-order existential quantifier and has a total of $n$ blocks of quantifiers.  We refer to a second-order theory as a $\forall_n$ theory if it has an axiomatization using sentences that are at most $\forall_n$, and likewise for $\exists_n$ theories.

\subsection{Descriptive complexity of second-order theories}
Morley's proof of his result in \cite{Morley1970} relies heavily on the fact that $\op{Mod}_S(\sigma)$ is Borel when $\sigma$ is an $L_{\omega_1, \omega}$ sentence.  In the context of a sentence of second-order logic the sets $\op{Mod}_S(\sigma)$ are no longer necessarily Borel.  Instead, the models of a single second-order sentence will form a projective set, with the complexity of the set being determined by the quantifier complexity of the sentence.

\begin{lem}\label{lem:ProjectiveComplexityOfSecondOrder}
Let $S$ be a countable signature, and let $\sigma$ be a second-order $S$-sentence.  Then $\op{Mod}_S(\sigma)$ is projective.  More specifically, if $\sigma$ is a $\forall_n$ formula then $\op{Mod}_S(\sigma)$ is a $\bPi_n^1$ set, and if $\sigma$ is an $\exists_n$ formula then $\op{Mod}_S(\sigma)$ is a $\bSigma_n^1$ set.
\end{lem}
\begin{proof}
The proof is by induction on the complexity of $\sigma$.  It is well-known (see e.g. \cite[Lemma 11.3.3]{Gao}) that if $\sigma$ is first-order then $\op{Mod}_S(\sigma)$ is Borel.  Conjunctions correspond to intersections and negations to complements, so the case of interest is the second-order existential quantifier.

Suppose that $\mathcal{M}$ is an $S$-structure, and $X$ is a second-order variable, and that $\mathcal{M}\models (\exists X) \varphi(X)$. Then there is an $S\cup\{X\}$-structure $\mathcal{M}'$ such that $\mathcal{M}'\models \varphi$ and its reduct to $S$ is $\mathcal{M}$, i.e. $\mathcal{M}'\restr S=\mathcal{M}$.  It is easily verified that the projection $f:\op{Mod}_{S \cup \{X\}} \to \op{Mod}_S$ given by $f(\mathcal{M})=\mathcal{M}\restr S$ is continuous.  We have:
\begin{align*}
    \op{Mod}_S(\exists X\varphi) &= \{\mathcal{M}\in \op{Mod}_S : \mathcal{M} \models (\exists X) \varphi(X)\} \\
    &=\{\mathcal{M} \in \op{Mod}_S : (\exists \mathcal{N}\in \op{Mod}_{S \cup \{X\}}) \mathcal{N}\models \varphi(X) \text{ and } \mathcal{N}\restr S=\mathcal{M}\} \\
    &=\{\mathcal{M}\in \op{Mod}_S : (\exists \mathcal{N}\in \op{Mod}_{S\cup\{X\}}(\varphi)) f(\mathcal{N})=\mathcal{M}\} \\
    &= f[\op{Mod}_{S\cup\{X\}}(\varphi)]
\end{align*}
Thus the second-order existential quantifier corresponds to projection, completing the proof.
\end{proof}

If $T$ is a second-order theory in a countable signature, then $T$ consists of at most countably many second-order sentences, say $T = \{\sigma_n : n < \omega\}$, and then $\op{Mod}_S(T) = \bigcap_{n < \omega}\op{Mod}_S(\sigma)$.  It is easy to see that the relation of isomorphism of countable structures is a $\bSigma_1^1$ relation on $\op{Mod}_S$.  Combining these observations with Lemma \ref{lem:ProjectiveComplexityOfSecondOrder} and the fact that the projective classes $\bPi^1_n$ and $\bSigma^1_n$ are closed under finite union (see \cite[Proposition 37.1]{Kechris}) we obtain the following result for second-order theories:

\begin{prop}\label{prop:SecondOrderComplexity}
If $T$ is a second-order theory of bounded quantifier complexity, then $\cong_T$ is a projective equivalence relation.  More specifically, for $n>1$, if $T$ is a $\forall_n$ theory then $\cong_T$ is a $\bSigma_n^1$ relation, and if $T$ is an $\exists_n$ theory then $\cong_T$ is a $\bPi_{n}^1$ relation. If $T$ is an existential second-order sentence (i.e. $\exists_1$), then $\cong_T$ is $\bDelta^1_2$.
\end{prop}

If we wish to consider second-order theories of unbounded quantifier complexity then it is useful to introduce the following definition that generalizes the projective hierarchy:

\begin{defi}
The collection of \textit{$\sigma$-projective sets} is the smallest $\sigma$-algebra containing the
open subsets (of $\R$) and closed under projections.
\end{defi}

\begin{lem}
Let $S$ be a countable signature, and let $T$ be a second-order $S$-theory.  Then $\op{Mod}_S(T)$ is $\sigma$-projective; it is in fact a countable intersection of projective sets.
\end{lem}
\begin{proof}
$\op{Mod}_S(T) = \bigcap_{\sigma \in T}\op{Mod}_S(\sigma)$, so this follows directly from Lemma \ref{lem:ProjectiveComplexityOfSecondOrder}.
\end{proof}

\section{Morley's Theorem fails consistently for second-order logic}\label{sec:ConsistentFailure}
Our strategy for showing the consistent failure of Morley's Theorem for second-order logic is to force over $L$ to add $\aleph_2$ Cohen reals, and then force over the resulting model to add $\aleph_3$ random reals.  After that forcing, we have $2^{\aleph_0} = \aleph_3$.  In this final model $L[G][H]$ of set theory there are exactly $\aleph_2$ Cohen reals over $L$ (random reals don't add Cohen reals: \cite[Section 7.2]{BJ95}), which will enable us to construct a second-order sentence with exactly $\aleph_2$ isomorphism classes of countable models.

\begin{lem}\label{lem:DST1}
Let $S = \{+, \cdot, <, 0, 1, R\}$, where $R$ is a unary relation symbol.  Suppose that $A \subseteq 2^\omega$ and $A$ is $\mathit{\Sigma^1_n}$ (respectively, $\mathit{\Pi^1_n}$) for some $n \geq 2$.  Then there is an $\exists_n$ (respectively, $\forall_n$) $S$-sentence such that every model of $\sigma$ is isomorphic to $(\omega, +, \cdot, <, 0, 1, R)$ for some $R \in A$, and moreover $(\omega, +, \cdot, <, 0, 1, R) \cong \\(\omega, +, \cdot, <, 0, 1, R')$ if and only if $R = R'$.  In particular, the number of isomorphism classes of models of $\sigma$ is $\abs{A}$.
\end{lem}
\begin{proof}
We give the proof for $\mathit{\Sigma^1_n}$; the proof for $\mathit{\Pi^1_n}$ is similar.

Let $PA_{II}$ be second-order Peano arithmetic, which can be expressed as a $\forall_1$ sentence of $S \setminus \{R\}$.  It is well-known that $PA_{II}$ is categorical, and (up to isomorphism) the only $S \setminus \{R\}$ structure satisfying $PA_{II}$ is $(\omega, +, \cdot, <, 0, 1)$.  Since $A$ is $\mathit{\Sigma^1_n}$, there is an $\exists_n$ formula $\phi(X)$, in one second-order variable $X$, such that for every $a \in 2^\omega$, $a \in A$ if and only if $(\omega, +, \cdot, < 0, 1) \models \phi(a)$ (see \cite[8B.15]{Mos80}; the second-order arithmetic there is described using Henkin semantics, but for the purposes of this argument the difference is only notational).

Let $\sigma$ be $PA_{II} \wedge \phi(R)$; $\sigma$ is $\exists_n$ because $PA_{II}$ is $\forall_1$ and $\phi$ is $\exists_n$ with $n \geq 2$.  As noted above, every model of $\sigma$ is isomorphic to one of the form $(\omega, +, \cdot, <, 0, 1, R)$ for some $R \in A$.  Finally, $(\omega, <)$ has no non-trivial automorphisms, so the only possible isomorphism from $(\omega, +, \cdot, <, 0, 1, R)$ to $(\omega, +, \cdot, <, 0, 1, R')$ is the identity map.
\end{proof}

\begin{thm}\label{mainthm}
It is consistent with $\ZFC$ that there exists a second-order sentence with exactly $\aleph_2$ non-isomorphic countable models while the continuum is $\aleph_3$.
\end{thm}
\begin{proof}
Force over $L$ to add $\aleph_2$ Cohen reals, and then force over the resulting model to add $\aleph_3$ random reals.  Let $C$ be the set of reals in this model that are Cohen over $L$.  Since adding random reals does not add reals that are Cohen over $L$ (see \cite[Section 7.2]{BJ95}), we have $\abs{C} = \aleph_2$, while $2^{\aleph_0} = \aleph_3$.  It is a folklore result that $C$ is $\mathit{\Pi^1_2}$.  See, for instance, \cite{HamkinsMO}, which gives an explicit $\forall_2$ definition of $C$ based on the fact that the set of reals in $L$ is $\mathit{\Sigma^1_2}$ (for which see \cite[8F.7]{Mos80}), and hence shows that $C$ is $\mathit{\Pi^1_2}$ by \cite[8B.15]{Mos80}.  Thus Lemma \ref{lem:DST1} provides a $\forall_2$ sentence with exactly $\aleph_2$ models (all of which are countable).
\end{proof}

Our original Cohen-random proof of Theorem \ref{mainthm} was fatally flawed. We thank the referee for spotting this and suggesting a proof along the lines we have done above.

\section{Morley's Theorem is consistently true for second-order logic}\label{sec:Consistency}
In this section we sketch that, modulo a supercompact cardinal, the extension of Morley's result to second-order logic is consistently true.  This completes the proof of the undecidability of the second-order version of Morley's Theorem.  We will obtain finer results, with smaller large cardinals, in later sections.  In fact, the consistency of the Absolute Morley Theorem follows directly from the following result of Foreman and Magidor.

\begin{thmr}[{\cite{FM95}}]\label{thm:ForemanMagidor}
If it is consistent that there is a supercompact cardinal, then it is consistent that $\neg\CH$ and
every equivalence relation on $\R$ that is in $L(\R)$ has $\leq \aleph_1$ or a perfect set
of inequivalent elements. In particular, this holds for $\sigma$-projective equivalence relations.
\end{thmr} 

Foreman and Magidor do not specifically state Theorem \ref{thm:ForemanMagidor} in \cite{FM95}, but it is implicit in their work.  In this section we will sketch a proof of the $\sigma$-projective version of Theorem \ref{thm:ForemanMagidor}.  Since our intended
audience includes model theorists who may not be familiar with large cardinals or inner models of set theory, we will need to briefly explicate ``supercompact'', ``$L(\R)$'', and other notions . But first, for the set theorist reader, we
should mention that with the development of Inner Model Theory since \cite{FM95}, it is now clear that
the supercompact cardinal above may be reduced to the assumption that there is a proper class of Woodin cardinals (defined in Section \ref{sec:Woodin}), or even a sequence of infinitely many Woodin cardinals with a measurable cardinal above all of them.  We will consider and prove various extensions and improvements of Theorem \ref{thm:ForemanMagidor} in Section \ref{sec:Woodin} below.

In Section \ref{sec:Preliminaries} we showed that when $\sigma$ is a second-order sentence, $\op{Mod}(\sigma)$ is projective; it is not generally Borel.  If $T$ is a second-order theory, then $\op{Mod}(T)$ is the countable intersection of projective sets, but these projective sets may have unbounded complexity, so $\op{Mod}(T)$ is $\sigma$-projective but may not be projective.  Since some results from the literature are set in the broader context of sets in $L(\mathbb{R})$, we remind the reader of that setting.

\begin{defi}
$L(\R)$, the collection of all sets constructible from $\R$, is the smallest inner model of $V$ containing $\R$ as a member. It is defined in analogy to $L$ and $L[\R]$ by:

\centering{$L_0(\R)=Trc(\R)$ the transitive closure of $\R$}
\\

\centering{$L_{\alpha+1}(\R)=Def(L_{\alpha}(\R))$}
\\

\centering{$L_\delta(\R)=\bigcup_{\alpha<\delta}L_{\alpha}(\R)$ for limit $\delta>0$}
\\

\centering{$L(\R)=\bigcup_{\alpha\in ON}L_{\alpha}(\R)$}

\end{defi}

We are actually interested in $L(\R)\cap \mathcal{P}(\R)$, but will usually just write $L(\R)$. Notice that
each element $C$ of $L(\R)$ can be defined by a formula of the language of set theory with a real
and finitely many ordinals as parameters. We speak of the real together with the ordinals as a
“code” of $C$. We will speak more about coding later.
In particular, by coding, we can easily see that $\sigma$-projective sets are in $L(\R)$. In fact,

\begin{propr}[{\cite[Proposition 3.8]{AMS21}}] The collection of $\sigma$-projective sets of reals is precisely $L_{\omega_1}(\R) \cap \mathcal{P}(\R)$.
\end{propr}

Weaker large cardinal hypotheses suffice for our $\sigma$-projective
purposes, as we shall see in Sections \ref{sec:Extension} and \ref{sec:Woodin}, and these are sufficient for our applications to second-order logic.  We do not know of any interesting logics that lead us beyond the $\sigma$-projective sets toward $L(\R)$, but any logic we
can define by using a formula of set theory with finitely many real and ordinal parameters is fair
game.
The model of set theory that Foreman and Magidor use is actually one familiar to set theorists. It
is the usual model for the \emph{Proper Forcing Axiom}.

\begin{defi}
    A poset $\mathbb{P}$ is \emph{proper} if and only if for all uncountable cardinals $\kappa$ and all stationary $S \subseteq [\kappa]^{\aleph_0}$, $1\Vdash_{\mathbb{P}} \mbox{``$S$ is stationary''}$.
\end{defi}

\begin{defi}
A cardinal $\kappa$ is \textit{supercompact} if for every cardinal $\lambda\geq \kappa$ there exists an
elementary embedding $j_\lambda$ of $V$ into an inner model $M$ (i.e. a proper class model of ZFC included in $V$) with \textit{critical
point} $\kappa$ (least ordinal $j_\lambda$ moves) and $\lambda < j_\lambda(\kappa)$, such that $M^\lambda$ is included in $M$.
\end{defi}
The Proper Forcing Axiom ($\PFA$) is like $\MA_{\aleph_1}$, but instead of meeting $\aleph_1$ dense
sets for countable chain condition partial orders, one does this for the wider class of proper partial orders. The
usual model for $\PFA$ is obtained by iterating proper posets of size less than a supercompact cardinal
$\kappa$ via countable support iteration $\kappa$ many times, and then using a reflection argument to argue that all collections of $\aleph_1$ dense sets in any proper partial order have been dealt with. For details, see e.g. \cite{Je03}. It is known that $\PFA$ implies $2^{\aleph_0} = \aleph_2$ (\cite[Theorem 1.8]{Velickovic}; see also \cite{Todorcevic}).

We will briefly sketch the ideas of the Foreman-Magidor proof in the special case we are
interested in, namely $\sigma$-projective equivalence relations. A detailed analysis of the proof will appear in Section \ref{sec:Woodin}, where we find weaker
hypotheses that still enable us to compute the possible numbers of equivalence classes. The
idea of the Foreman-Magidor proof is to note that a code for such a set involving a real plus finitely many countable ordinals appears at some initial stage of the iteration; if such an
equivalence relation is \emph{thin} -- i.e., does not have perfectly many equivalence classes -- then they show that the interpretation of the code when the code
appears is also thin (``downwards generic absoluteness”), and that the rest of the proper forcing
adds no new equivalence classes to thin equivalence relations (``upwards generic
absoluteness”). When the code appears at some initial stage, we may without loss of generality
assume $\CH$ holds at that stage, since it holds cofinally often in the iteration. At that stage then,
its interpretation has no more than $\aleph_1$ equivalence classes. By upwards generic
absoluteness, its interpretation at the final stage, which is just the equivalence relation we
started with, then has no more than $\aleph_1$ equivalence classes. It is convenient that we can
refer to a well-known model of set theory, but the Foreman-Magidor proof does not actually use
most of the properties of that model.  It is actually sufficient to alternately blow up the continuum, properly, e.g. with Cohen reals, and then countably closedly collapse it down to $\aleph_1$, in a countable support iteration.  For further discussion, see Section \ref{sec:Woodin}.

In Section \ref{sec:Woodin}, with the benefit of research after \cite{FM95}, we will prove results sharper than Theorem \ref{thm:ForemanMagidor}.  We need to talk about codes. Think for example of the pair $\langle m,n\rangle$ as the real $2^{m+1} \times 3^{n+1}$. This real can be
thought of as coding the open interval $(m,n)$, which of course has different extensions in different
models. $\sigma$-projective sets are coded by a real, just as Borel sets are. The reader is
probably familiar with the idea of coding a Borel set: one first lists the pairs of rational numbers
in some recursive way, then countable sequences of such pairs to code open sets, then
countable sequences of these, as well as countable intersections of complements of these, then
continue recursively to list countable sequences of what has gone before, etc. For a careful exposition of coding, see \cite{Kun11}.  Sets in $L(\R)$ are
coded by a real plus finitely many ordinals. This makes things more complicated, since those
ordinals may be bigger than the large cardinal that we are collapsing to $\aleph_2 = 2^{\aleph_0}$. Since we
don’t currently have any examples of interesting equivalence relations that are in $L(\R)$ but are not definable
from a real, here we have just sketched the Foreman-Magidor proof for equivalence relations that are in $L(\R)$ and definable from
a real plus finitely many \emph{countable} ordinals, i.e. the $\sigma$-projective equivalence relations. We will prove -- indeed improve -- their theorem in Section \ref{sec:Woodin}.

\section{Finer Analysis: Descriptive Set Theory}\label{sec:Extension}
In the previous section we noted that from large cardinals we can obtain the consistency of the Absolute Morley Theorem for second-order logic.  In this section we start examining how weakenings of the $L(\mathbb{R})$ results of Foreman-Magidor \cite{FM95} can be achieved via weaker assumptions.  We achieve this by cobbling together results from the literature.  In Section \ref{sec:Woodin} we employ the machinery of Inner Model Theory to achieve more precise conclusions. 

Most of this section concerns Descriptive Set Theory, and as such when we discuss issues of definability we mean ``definable" in the sense used in Descriptive Set Theory.  In particular, typical definability notions we consider are ``Borel'', ``analytic'', ``projective'', or being in $L(\mathbb{R})$. 

Morley and Burgess use that projections of Borel sets are analytic.  To extend their results to second-order logic we can use that countable
intersections of projective sets are $\sigma$-projective. Morley quotes the classical theorem of Descriptive Set Theory that
analytic sets have the perfect set property; we want to make similar assertions for $\sigma$-projective sets. However such assertions are no longer theorems of
$\ZFC$: Gödel showed that under $\VL$, there is an uncountable co-analytic ($\bPi^1_1$) set which
does not include a perfect set. See e.g. \cite[Theorem 13.12]{Ka08} for a proof. Co-analytic sets do, however, have
cardinality $\leq \aleph_1$ or $2^{\aleph_0}$.

One way of extending the perfect set property to more complicated sets of reals is via large cardinals, which imply restrictions of the Axiom of Determinacy to various classes of definable sets of reals, which in
turn imply that if such sets are uncountable, they include perfect sets. We refer the reader to
Section 27 of \cite{Ka08} for an introduction to determinacy and the formal definitions of games and
winning strategies. In the standard Descriptive Set Theory abuse of notation, we think of $\R$ or $[0, 1]$ as
$\omega^\omega$. This is harmless since the descriptive structures of these spaces are the same; that is, these spaces, as well as the Cantor space, are Borel-isomorphic. 

\begin{defi}
Let $\mathcal{C}$ be a collection of subsets of $\R$ (we are interested in $\mathcal{C}$’s which are composed of sets which are “definable” in some sense). $\AD_\mathcal{C}$ is the assertion that given any $C \in \mathcal{C}$, in any game where players alternately pick natural numbers, with Player I trying to get the resulting infinite sequence to be in $C$ and Player II trying to prevent that, one of the players has a winning strategy. \emph{Projective determinacy} ($PD$) is the assertion $\AD_{\mathcal{P}}$ where $\mathcal{P}$ is the collection of all projective sets. For specific levels of the projective hierarchy, instead of e.g. $\AD_{\bSigma^1_{17}}$, it is common to write $Det(\bSigma^1_{17})$.
\end{defi}

Letting $\mathcal{B}$ stand for the collection of all Borel sets, $\AD_\mathcal{B}$ is a theorem of $\ZFC$ (\cite{Mar75}, see also \cite[Section 20]{Kechris}); $\AD$ for progressively larger definable classes follows from progressively larger cardinals. $L(\R)$ cannot satisfy the Axiom of Choice if we suppose that it satisfies the Axiom of Determinacy, but there is much to be said for assuming $\AD_{L(\R)}$
within a $\ZFC$ environment, since it follows from large cardinals and imposes a pleasant
regularity on the sets of reals constructible from $\mathbb{R}$.  See \cite[Section 32]{Ka08}. For example,

\begin{thm}\label{thm:MSW} 
If there is a supercompact cardinal, then $\AD_{L(\R)}$ holds. 
\end{thm}

The hypothesis of Theorem \ref{thm:MSW} has by now been considerably weakened, replacing the supercompact cardinal with a measurable cardinal above infinitely many Woodin cardinals \cite[Corollary to Main Theorem and Woodin Theorem]{Martin1989}. The definition of a Woodin cardinal can be found in Section \ref{sec:Woodin}. For details, see e.g. \cite{Neeman2010}. This hypothesis is considerably weaker than the existence of a supercompact cardinal. For sufficiently closed classes
$\mathcal{C}$, $\AD_\mathcal{C}$ implies every uncountable member of $\mathcal{C}$ includes a perfect set \cite[Theorem 4.1]{Davis}.

We want to build a model of set theory in which the isomorphism relation for countable models of a second-order theory has either at most $\aleph_1$ classes or has a perfect set of non-isomorphic models.  We don’t know of any results that deal specifically with this question, but there is a long line of research extending Silver \cite{Sil80} which counts the number of equivalence classes of a definable equivalence relation among sets of
reals.

\subsection{Thin equivalence relations}\label{sec:Thin}
A great deal is known about equivalence relations on Polish spaces that do, or do not, have a perfect set of pairwise inequivalent elements.  Here we remind the reader of some of the standard terminology and some relevant results from the literature.  We also prove that Morley's Theorem extends to sentences of universal second-order logic.  

Although results about determinacy and the perfect set property are motivating, there is no link between a projective class having the perfect set property and whether an equivalence relation of that complexity has perfectly many equivalence classes. In fact, $\bSigma^1_2$ has the perfect set property because each $\bSigma^1_2$ set is the union of $\aleph_1$ Borel sets (a classical result), but it is undecidable, modulo large cardinals, whether
or not a $\bSigma^1_2$ equivalence relation must have $\leq \aleph_1$ or perfectly many equivalence classes.

Recall from Section \ref{sec:Preliminaries} that an equivalence relation $E$ on a Polish space $X$ is \textit{thin} if there is no perfect set of pairwise $E$-inequivalent elements of $X$.  The theorem of Burgess \cite{Burgess} mentioned in the introduction can thus be stated as follows:

\begin{thm}[Burgess]\label{thm:Burgess}
A $\bSigma^1_1$ thin equivalence relation on a Polish space $X$ has at most $\aleph_1$ equivalence classes.
\end{thm}

As an immediate consequence, we obtain the following generalization of Morley's result.

\begin{crl}\label{crl:UniversalSecondOrder}
Let $S$ be a countable signature, and let $\sigma$ be a universal second-order $S$ sentence (that is, a sentence for which the second-order quantifiers in the prenex form of $\sigma$ are all universal).  Either there is a perfect set of non-isomorphic models of $\sigma$, or there are at most $\aleph_1$ non-isomorphic models of $\sigma$.
\end{crl}
\begin{proof}
By Proposition \ref{prop:SecondOrderComplexity} the relation $\cong_\sigma$ is $\bSigma_1^1$, and hence Burgess' theorem applies.
\end{proof}

Schlicht's Example \ref{ex:Schlicht} below establishes limits to what we can prove in ZFC plus large cardinals about the number of equivalence classes of thin projective equivalence relations. We shall see below, however, that large cardinals do enable us to extend ZFC results, since they imply some determinacy. See e.g. Theorem \ref{pdcohen} and Corollary \ref{supercompact}. Of course generic extensions of large cardinal models give us even more power: viz. \cite{FM95}.

\begin{ex}
Let $S = \{<\}$, and let $\sigma$ be the second-order sentence expressing that $<$ is a well-order.  That is, $\sigma$ is the conjunction of the (first-order) axioms for being a linear order with the second-order statement 
\[(\forall^1 A)~(((\exists z)~A(z)) \to (\exists x)~(A(x) \wedge (\forall y)(A(y) \to x \leq y))).\]
Up to isomorphism, the countable models of $\sigma$ are exactly the countable ordinals, so $\sigma$ has $\aleph_1$ pairwise non-isomorphic models, regardless of the size of the continuum.  This example shows that Vaught's Conjecture is consistently false  (and hence undecidable, since it is true assuming the continuum hypothesis) for second-order logic, and even for single sentences of universal second-order logic.

For this example, $\sigma$ is a universal second-order statement, so $\cong_\sigma$ is $\bSigma_1^1$ by Proposition \ref{prop:SecondOrderComplexity}.  This example shows that Silver's Dichotomy for $\bPi_1^1$ equivalence relations (see \cite[Theorem 5.3.5]{Gao}) does not extend to $\bSigma_1^1$ relations.  Since isomorphism of countable structures is already $\bSigma_1^1$ even without restricting to models of a specific sentence, we see that there is little hope for directly applying Silver's Dichotomy to make progress on Vaught's conjecture, even for first-order theories.
\end{ex}

The situation for existential second-order sentences is no better than it is for universal sentences, as shown by the next example, which is due to Kunen (see \cite[1.4.3]{Steel}).

\begin{ex}
Let $S = \{<\}$, and let $\sigma$ be the second-order sentence obtained by taking the conjunction of the (first-order) axioms for $<$ being a linear order with the second-order statement expressing that for any two elements $a$ and $b$ there is an automorphism of $<$ sending $a$ to $b$.  Kunen shows that every model of $\sigma$ is of the form $Z^\alpha$ or $Z^\alpha \cdot \eta$ for some countable ordinal $\alpha$, where $Z = \{0\}$, the product order is lexicographic, and $\eta$ is the ordertype of $\mathbb{Q}$.  In particular, $\sigma$ has exactly $\aleph_1$ countable models (and Steel \cite{Steel} also points out that $\sigma$ does not have a perfect set of non-isomorphic countable models).

For this example, $\sigma$ is an existential second-order statement, so $\cong_\sigma$ is $\bDelta^1_2$ by Proposition \ref{prop:SecondOrderComplexity}.
\end{ex}

The third example we wish to consider is due to Schlicht \cite{Schlicht}.  Some background is required before we can state the example.  Recall that a \emph{prewellordering} on a set $X$ is a relation $\leq$ on X which is reflexive, transitive, connected (i.e. for any $x,y \in X$, $x \leq y$ or $y \leq x$), and every nonempty subset of $X$ has a least element.  There is a natural equivalence relation associated with a prewellorder: $x \sim y$ if and only if $x \leq y$ and $y \leq x$.  It is interesting to note that, as pointed out in \cite[Remark 5.25]{Schlicht}, Projective Determinacy ($PD$) and the Axiom of Dependent Choice ($DC$) together imply that every thin projective equivalence relation is induced by a projective prewellorder.

For $n \geq 1$, the $n$th \emph{projective ordinal} $\delta^1_n$ is the supremum of lengths of $\bDelta^1_n$
prewellorders.  The statement that ``$x^\#$ exists" is equivalent to there being a non-trivial elementary embedding
from $L[x]$ to $L[x]$ - see e.g. \cite[Section 9]{Ka08}.  If there is a measurable cardinal, or even just a
cardinal $\kappa$ such that $\kappa \to (\omega)^{<\omega}_{2}$, then $x^\#$ exists for every real $x$. $PD + DC$ also implies this conclusion; indeed, over $ZF+DC$, analytic determinacy is equivalent to the existence of $x^\#$ for every real $x$ - see \cite[Theorem 33.19]{Je03}. $PD$ is consistent with $AC$ and follows from the existence of infinitely many Woodin cardinals \cite[Corollary to Main Theorem]{Martin1989}.

We are now ready to describe Schlicht's example.  For more details, see \cite[Example 4.18]{Schlicht}.
\begin{ex}\label{ex:Schlicht}
Assume $x^\#$ exists for every real $x$. Let $\langle \iota^x_\alpha: \alpha \in ORD\rangle$ enumerate
the $x$-indiscernibles and define $u^x_2 = \iota^x_{\omega_1+1}$. The prewellorder defined by
$$x \leq y \iff u^x_2 \leq u^y_2$$ 
is $\bDelta^1_3$ and has length $\delta^1_2$. On the other hand, Kunen-Martin  (unpublished, but see \cite[2G.2]{Mos80}) proved from $ZF+DC$ that $\delta^1_1 = \aleph_1$, while assuming $ZF+DC$, Martin proved that $\delta^1_2 \leq \aleph_2$ (see these results in \cite[p.162]{Mos72}). But, also assuming $ZF+DC+PD$, Kechris and Moschovakis proved $\delta^1_1 < \delta^1_2$ (the results of Martin and Kechris-Moschovakis are found in \cite[Theorem 9.1]{KechrisMoschovakis}). Thus,
assuming $ZF+DC+PD$, we get a thin $\bDelta^1_3$ equivalence relation with exactly $\aleph_2$ equivalence classes.

This example shows that the result of Foreman and Magidor (used in Section \ref{sec:Consistency}) in the usual PFA model does not hold in ZFC plus a sufficiently large cardinal.  Indeed, suppose for example that we have a supercompact cardinal. Add $\aleph_3$ Cohen reals. The supercompact is still supercompact so $PD$ holds and we have sharps, so Schlicht's $\bDelta^1_3$ example exists, and has $\aleph_2$ equivalence classes, and is thin because $\aleph_2 < \aleph_3 = 2^{\aleph_0}$.
\end{ex}

It is unclear if Schlicht's equivalence relation is of the form $\cong_T$ for a second-order theory $T$, so we can still ask:

\begin{prob}
Prove that a supercompact cardinal does not imply (with $2^{\aleph_0}  > \aleph_2$) Absolute Morley for second-order theories.
\end{prob}

In between Burgess’ dichotomy for $\bSigma_1^1$ equivalence relations and Schlicht's $\bDelta^1_3$ example, one might wonder about
$\bSigma^1_2$ equivalence relations, especially in view of the Shoenfield Absoluteness Theorem. The following result of Harrington-Sami is
encouraging:

\begin{lemr}[{\cite[Theorem 5]{HSa}}]\label{lem:HSa}
Assume $PD$. Thin $\bSigma^1_2$ equivalence relations are $\bDelta^1_2$.
\end{lemr}

However, Schlicht’s analysis of the number of classes of projective equivalence
relations specializes in the $\bSigma^1_2$ case to yield:

\begin{thm}\label{thm:Schlicht}
In $ZFC$ 
there is no upper bound below $2^{\aleph_0}$ for the number of equivalence classes of thin
$\bSigma^1_2$ equivalence relations.
\end{thm}

\begin{proof}
  Given a set $A$ of reals, consider the trivial equivalence relation $E(A)$ on $\mathbb{R}$ defined by:
\[x \sim y \text{ if and only if } x = y \text{ or } x \notin A \text{ or } y \notin A.\]

Then $E(A)$ has $\abs{A}+1$ equivalence classes. If $A$ is $\bPi^1_2$, then $E(A)$ is $\bSigma^1_2$. Now, in the proof of Theorem \ref{mainthm}, rather than adding $\aleph_2$ Cohen reals and $\aleph_3$ random reals, we  could have added $\kappa$ many Cohen reals and $\kappa^+$ many random reals, for any regular $\kappa$, to get $2^{\aleph_0}$  to be $\kappa^+$ while there is a thin $\bSigma^1_2$ equivalence relation with $\kappa$ many equivalence classes.
\end{proof}

Theorem \ref{thm:Schlicht} is actually a special case of Schlicht's Lemma 4.17 in \cite{Schlicht}, but that is not so easy to see, since his statement and proof involve Woodin cardinals, determinacy, and \emph{premice}. He uses Harrington forcing from \cite{Har77} instead of random forcing.

Thus, even for the relatively simple case of $\bSigma^1_2$ we must either assume large cardinals or, as in the Foreman-Magidor theorem, consider generic
extensions of the universe in order to get a Burgess-type result. What we can get from known results is:

\begin{thm}\label{pdcohen}
Assume $PD$.  Then thin $\bSigma^1_2$ equivalence relations have at most $\aleph_1$ equivalence classes.
\end{thm}
\begin{proof}
Harrington-Shelah \cite[Corollary 3]{HSh} showed the conclusion for $\bPi^1_2$, assuming there was a Cohen real over $L$.  They then applied Lemma \ref{lem:HSa}, getting $\bSigma^1_2$
thin equivalence relations are $\bPi^1_2$.  But $PD$ - indeed $\Det(\bPi^1_1)$ - implies there is a Cohen real over $L$ because it implies $\aleph_1^L$ $(= (2^{\aleph_0})^L)$ is countable, so there are only countably many dense sets to meet.  See more on Harrington-Shelah in Sections \ref{Suslin} and \ref{game}.
\end{proof}


We do not know if Theorem \ref{pdcohen} is already known. Kechris \cite[proof of Theorem 3.2]{Kechris1978} proved the conclusion from $\ZFC+\AD_{L(\R)}$. Either way, from Theorem \ref{thm:MSW} we have:

\begin{crl}\label{supercompact}
  If there is a supercompact cardinal, then thin $\bSigma^1_2$ equivalence relations have at most $\aleph_1$ equivalence classes.
\end{crl}

\begin{rmk}
Theorem \ref{pdcohen} and Corollary \ref{supercompact} do not contradict Theorem \ref{mainthm}: in a universe in which PD holds, our Cohen plus random forcing over $L$ only creates an equivalence relation with $\aleph_2^L$, i.e. countably many, equivalence classes.
\end{rmk} 

\begin{rmk}
There is another interesting way of looking at Theorem \ref{thm:Schlicht}.  Harrington \cite[Theorem B]{Har77} proved:
\end{rmk}
\begin{lem} There is a model of ZFC in which:
\begin{enumerate}
    \item $2^{\aleph_0}$ can be as large as one likes,
    \item Martin’s Axiom holds,
    \item Every set of reals of size less than $2^{\aleph_0}$ is $\bPi^1_2$.
\end{enumerate}
\end{lem}

In Harrington’s model, for each $\kappa < 2^{\aleph_0}$, there are thin $\bSigma^1_2$ equivalence relations with $\kappa$-many equivalence classes. Again, this does not contradict Theorem \ref{pdcohen} because Harrington’s model is constructed as a cardinal-preserving extension of a model in which $\aleph_1 = \aleph_1^L$, which contradicts PD.

\subsection{$\aleph_1$-Suslin and Co-$\aleph_1$-Suslin Equivalence Relations}\label{Suslin}
In Foreman-Magidor \cite{FM95} and other works involving determinacy, \emph{$\kappa$-Suslin} sets and their complements play a role. They also appear in the study of thin equivalence relations.  

\begin{defir}[{\cite[p. 148]{HSh}}]\mbox{}
\begin{enumerate}[label=(\alph*)]
  \item $T$ is a \emph{tree} on the set $Y$ if: $T\subseteq Y^{<\omega}$, and $(\eta \in T, \tau\subseteq \eta) \implies \tau \in T$.
  \item For $T$ a tree on $Y$, $[T]:= \st{f}{f: \omega \to Y, (\forall n)(\rstrct{f}{n} \in T)}$.
  \item For $T$ a tree on $\kappa \times X$, let
  \[p[T]:= \st{g}{g:\omega \to X, \mbox{ and for some } h:\omega \to \kappa, \pair{h}{g} \in [T]}.\]
  (Here we identify $\pair{h}{g}$ with the function $f:\omega \to \kappa\times X$ where $f(n) = \pair{h(n)}{g(n)}$).
  \item A binary relation $R$ on $\omega^\omega$ is \emph{$\kappa$-Suslin (via $T$)} if: $T$ is a tree on $\kappa\times(\omega^2)$ and $R = p[T]$; ($R$ is \emph{co-$\kappa$-Suslin} if: $\comp{R}$ is $\kappa$-Suslin (where $\comp{R} = $ complement of $R$)).
\end{enumerate}
\end{defir}

The following theorem follows easily from results in \cite{HSh}.

\begin{thm}
  $\MA_{\omega_1}$ implies thin co-$\aleph_1$-Suslin equivalence relations on $\R$ have $\leq\aleph_1$ equivalence classes. 
\end{thm}

\begin{proof}
  Since the tree $T$ for the complement of the co-$\aleph_1$-Suslin equivalence relation $E$ has cardinality $\aleph_1$, it follows that $L[T] \models \CH$, so ${\abs{2^{\omega}}}^{L[T]} = \aleph_1$, which $\MA_{\omega_1}$ implies is $<2^{\aleph_0}$. Then $\MA_{\omega_1}$ easily implies there is a real Cohen-generic over $L[T]$. That, according to \cite{HSh}\footnote{The first proof of this in \cite{HSh} has a gap, but the second is OK.}, implies $E$ remains an equivalence relation after a Cohen real is adjoined to $L[T]$. We now need another result from \cite{HSh}:

  \begin{lem}[{\cite[Theorem 1]{HSh}}]
    Suppose $E$ is a thin co-$\kappa$-Suslin relation via $T$, i.e.~$T$ is the tree for the complement of $E$. Assume $E$ is an equivalence relation after adding a Cohen real to $L[T]$. Then $E$ has at most $\kappa$ equivalence classes. \qedhere
  \end{lem}
 \end{proof}
  
  There is a connection between (co)-$\kappa$-Suslin and our previous discussion of $\bSigma^1_2$ and $\bPi^1_2$. First we quote:
  
   \begin{lemr}[{\cite[Paragraph 1.2 and Theorem 3.5]{Mar11}}]\label{lem:SuslinBorel}
    A relation $R$ is $\aleph_1$-Suslin if and only if it is the union of $\aleph_1$ Borel sets.
  \end{lemr}

The relationship between $\aleph_1$-Suslin and $\bSigma^1_2$ is quite interesting. It is a classical result that

\begin{lem}
Every $\bSigma^1_2$ set is the union of $\aleph_1$ Borel sets, i.e. is $\aleph_1$-Suslin.

\end{lem} 
Thus $\bPi^1_2$ sets are co-$\aleph_1$-Suslin, so we have

\begin{crl}
$MA_{\omega_1}$ implies thin $\bPi^1_2$ equivalence relations have $\leq \aleph_1$ equivalence classes.

\end{crl}

\begin{crl}\label{cor:MAExistential}
  $\MA_{\omega_1}$ implies an existential second-order theory has either a perfect set of pairwise non-isomorphic countable models or has $\leq \aleph_1$ of them.
\end{crl}

  On the other hand, Martin and Solovay \cite[p. 166]{MS} prove:

  \begin{lem}\label{lem:MS}
    Assume $\MA_{\omega_1}$. There is a $t \subseteq \omega$ such that $\omega_1 = \omega_1^{L[t]}$ if and only if the union of $\aleph_1$ Borel sets is $\bSigma^1_2$.
  \end{lem}
  
 Interestingly, their use of almost disjoint coding to make an $\aleph_1$-Suslin set $\bSigma^1_2$ is also found in Harrington \cite{Har77}, which Schlicht used instead of random forcing to prove Theorem \ref{thm:Schlicht}.

\begin{prob}
Is it consistent that CH fails but that thin $\aleph_1$-Suslin equivalence relations have $\leq \aleph_1$ equivalence classes?  
\end{prob}

One is tempted to apply Theorem \ref{pdcohen} plus Lemma \ref{lem:MS}, but their hypotheses are not mutually consistent.

\section{Other logics}\label{sec:OtherLogics}
In this section we consider the applicability of our results to logics other than second-order logic.
\subsection{Game quantifiers}\label{game}
Moschovakis \cite{Mos72, Mos80} and others, e.g.~\cite{Vaananen}, have considered well-ordered quantifier strings of length $\omega$. The question of the number of equivalence classes for theories involving these fits in nicely with topics we have been discussing. In particular, we have the closed and open game quantifiers discussed in \cite{Mos72}:

\begin{defi}
The expression
$$(\forall x_0\exists y_0 \forall x_1 \exists y_1...) \bigwedge_{n<\w}\varphi_n(x_0,y_0,...,x_n,y_n)$$
is a \textit{closed game quantifier sentence of length $\w$}. The truth value of this game expression in a model $\mathcal{M}$ is equivalent to the existence of a winning strategy for Player II in the following game:\\
\[ \begin{array}{c|ccccc} \mathrm{I} & a_0 & & a_1 & &\hdots \\ \hline
    \mathrm{II} & & b_0 & & b_1 & \hdots 
   \end{array} \]
   
   Here $a_0, b_0, a_1, b_1, \dots$ are elements of $M$ and $b_0, b_1, \dots$ are chosen such that \[ \mathcal{M}\models \varphi_0(a_0,b_0), \] \[ \mathcal{M}\models \varphi_1(a_0,b_0,a_1,b_1), \] etc. (otherwise Player II loses).
   
\end{defi}

The open game quantifier interchanges the universal and existential quantifiers and exchanges the infinite conjunction for an infinite disjunction. Thus, the negation of a closed game quantifier sentence is an open game quantifier sentence. Explicitly, an \textit{open game quantifier sentence of length} $\omega$ has the following form:

$$(\exists x_0\forall y_0 \exists x_1 \forall y_1...) \bigvee_{n<\w}\varphi_n(x_0,y_0,...,x_n,y_n).$$

The truth value of an open game quantifier sentence in a model $\mathcal{M}$ is equivalent to the existence of at least one play of the game in which Player I wins; equivalently, Player II does not have a winning strategy.  Note that, unless there are enough determinacy assumptions, this does not necessarily imply that Player I does have a winning strategy. 

\begin{defir}[\cite{Mos80}]
  A \emph{pointclass} $\Gamma$ is a collection of sets such that each element of $\Gamma$ is a subset of some finite product $X$ of Polish spaces. A pointclass is \emph{adequate} if it contains all recursive pointsets and is closed under recursive substitution, $\wedge$, $\vee$, and bounded existential and universal quantification. For simplicity, all our pointclasses are closed under Polish continuous preimages.
\end{defir}

For example, each level $\bSigma^1_n$ or $\bPi^1_n$, $n<\omega$, of the projective hierarchy is adequate.

Now we introduce the game operator $\Game$ following Moschovakis. The reader is referred to \cite{Mos80} for a thorough exposition on the subject. 

\begin{defi}

\begin{enumerate}
    \item Given a set $P\subseteq \script{X}\times \w^\w$, we define the set $\opgm P$ as follows: $x\in \opgm P \iff$ Player I wins the game $\{ \alpha : P(x,\alpha)\}$, i.e. the game with parameter $x$ where Player I plays each $\alpha(2n)$ and Player II plays each $\alpha(2n+1)$, and Player I wins if and only if $P(x,\alpha)$. 
    \item Given a pointclass $\Gamma$, we define $\opgm \Gamma =\{\opgm P : P\subseteq \script{X}\times \w^\w\wedge P\in\Gamma \}$.
\end{enumerate}

\end{defi}

The following theorem relates the game operator $\opgm$  and the open game quantifier: 

\begin{thm}
Given a language, fix a theory $T$ and a sequence of formulas $\langle \varphi_n : n<\w\rangle$ for which there is an adequate pointclass $\Gamma$ such that for any $\mathcal{M}\in\Mod(T)$ and $a_1, a_2,... \in M$, \[ \Mod(T\wedge\bigvee_{n<\w}\varphi_n(a_1,..., a_{2n}))\in{\Gamma} \] (as a set of reals coding the respective models). Then the formula $\psi$ defined by applying an open game quantifier to the sequence $\langle \varphi_n : n<\w\rangle$ satisfies \[ \Mod(\psi)\in\opgm\Gamma. \]
\end{thm}
\begin{proof}
For simplicity, we assume that our models are enumerated and that their elements are exactly the elements of $\w$. Define $P$ to be the following set:
\begin{align*}
    \{(\mathcal{M},&\alpha)\in \Mod(T)\times \w^\w : \mathcal{M}\models T \; \wedge \\ & (\exists n) (0<n<\w) \; \mathcal{M}\models \varphi_n (\alpha(0), \alpha(1),..., \alpha(2n-1))\}.
\end{align*}

Notice that for a fixed model $\mathcal{M}$ of $T$, the projection $\{\alpha : (\mathcal{M},\alpha)\in P\}$ corresponds to the set of plays which result in Player I winning the game corresponding to $\psi$. Thus $\opgm P$ is the set of models for which Player I wins the game and so $\opgm P=\Mod(\psi)$.

Finally, notice that $(\mathcal{M},\alpha)\in P$ if and only if
$$\mathcal{M}\models T \wedge \bigvee_{n<\w}\varphi_n(\alpha(0), \alpha(1),..., \alpha(2n-1))$$
That is $(\mathcal{M},\alpha)\in P$ if and only if $\mathcal{M}\in \Mod(T\wedge\bigvee_{n<\w}\varphi_n(\alpha(0),..., \alpha(2n-1)))$. Thus $\Mod(\psi)=\opgm P\in \opgm \Gamma$.
\end{proof}

\begin{crl}
 Let $\langle \varphi_n : n<\w\rangle$ be a sequence of second-order formulas of uniformly bounded complexity, i.e. there is an $m<\w$ such that every $\varphi_n$ has complexity at most $\bSigma^1_m$. If $\psi$ is the formula obtained by applying an open game quantifier to $\langle \varphi_n : n<\w\rangle$, then $\Mod(\psi)$ has complexity at most $\opgm \bSigma^1_m$. 
\end{crl}

We need the following notation from \cite{Mos80}: 

\begin{defi}
Suppose $\Gamma$ is a pointclass and $P\in\Gamma$, $P\subseteq \script{X}\times \w^\w$. We denote by $\exists^1P$ the set $\{x\in \script{X} : (\exists^1y)P(x,y)\}$ where 
``$\exists^1y$'' is short for ``$\exists^1y \in {^{\omega}{\omega}}$''. We define the pointclass $\exists^1\Gamma$ to be $\{\exists^1P : P\in\Gamma\}$. The definitions of $\forall^1P$ and $\forall^1\Gamma$ are analogous. 
\end{defi}

Using the notation in the previous definition, we can use the expression $\exists^1\Gamma\subseteq \Gamma$ to denote that $\Gamma$ is closed under second-order existential quantification, e.g. $\exists^1 \bSigma^1_2\subseteq \bSigma^1_2$.

Under further assumptions, classes of the form $\opgm \bSigma^1_m$ are well understood. Moschovakis \cite{Mos80} proves:

\begin{lemr}[{\cite[Theorem 6D.2(vi)]{Mos80}}]
Suppose $\Gamma$ is an adequate pointclass and $\Det(\Gamma)$ holds. If $\exists^1 \Gamma\subseteq \Gamma$, then $\opgm \Gamma = \forall^1\Gamma$. 
\end{lemr}

With further work, the following result is obtained:

\begin{lemr}[\cite{Mos80}]
  If $\PD$ holds, then
  \[\opgm\bSigma^0_1 = \bPi^1_1,\stp \opgm\bPi^1_1 = \bSigma^1_2,\stp \opgm\bSigma^1_2 = \bPi^1_3, \ldots\]
\end{lemr}

Thus, assuming $\PD$, the open and closed game quantifiers do not lead us out of the projective sets. Moreover, since the class of $\sigma$-projective sets is an adequate pointclass that is closed under second-order existential and universal quantifications by definition, we have:

\begin{crl}
If $\sigma$-projective determinacy holds, then $\opgm A$ is $\sigma$-projective for any $\sigma$-projective $A$.
\end{crl}

In particular, we can also obtain a $\ZFC$ result on isomorphism relations since Borel determinacy holds in $\ZFC$:

\begin{crl}
  If a sentence $\sigma$ consists of an open game quantifier followed by a sequence of $L_{\omega_1, \omega}$ formulas, then $\Mod(\sigma)$ is $\bPi^1_1$, so the associated equivalence relation is $\bSigma_1^1$.  Thus if the relation $\cong_\sigma$ is thin, there are $\leq\aleph_1$ equivalence classes.
\end{crl}

\begin{proof}
  In this situation $\Mod(\sigma)$ is obtained by applying $\opgm$ to a Borel set, and hence is $\bPi^1_1$; it follows that $\cong_\sigma$ is $\bSigma_1^1$ by Lemma \ref{lem:ProjectiveComplexityOfSecondOrder}.  The result on the number of equivalence classes then follows by Burgess' dichotomy theorem (Theorem \ref{thm:Burgess}).
\end{proof}

\begin{crl}\label{thm:UsualPFAGame}
  In the usual model for $\PFA$, if $\mathcal{L}$ is the closure of second-order logic under open and closed game quantifiers and $\sigma$ is an $\mathcal{L}$-sentence, then the associated isomorphism relation $\cong_\sigma$ has either a perfect set of pairwise non-isomorphic elements or $\leq \aleph_1$ of them.
\end{crl}

\begin{proof}
    The relation $\cong_\sigma$ is in $L(\mathbb{R})$, indeed it is $\sigma$-projective, so the result follows from Foreman-Magidor (Theorem \ref{thm:ForemanMagidor}).
\end{proof}

We can actually do better than this. See Theorem \ref{thm:WoodinsWithMeasurableAbove} below to see that we can get the conclusion of Corollary \ref{thm:UsualPFAGame} from a model which only assumes that there are infinitely many Woodin cardinals with a measurable cardinal above.

For one application of the closed game quantifier to a first-order formula, we get that the complexity of the equivalence relation is $\bDelta^1_2$, assuming $\Det(\bPi_1^1)$. We can now prove:

\begin{thm}\label{thm:measurablePi11}
  If $\Det(\bPi_1^1)$ then for a thin equivalence relation defined by a sequence of $L_{\omega_1,\omega}$ formulas preceded by a closed game quantifier, the number of equivalence classes is $\leq\aleph_1$.
\end{thm}

\begin{proof}
    By \cite[11.5, 27.13]{Ka08}, $\Det(\bPi^1_1)$ implies $\aleph_1^{L[a]}$ is countable for all reals $a$. By \cite{HSh} that implies $\bPi^1_2$ thin equivalence relations have $\leq \aleph_1$ equivalence classes.
\end{proof}

It is known that $\Det(\bPi^1_1)$ follows from the existence of a measurable cardinal.  This fact is commonly attributed to Martin \cite{Mar70}, though it does not appear explictly in that paper.  For the connection between $\Det(\bPi^1_1)$, measurable cardinals, and the results of \cite{Mar70}, see \cite[Exercises 33.11, 33.12]{Je03}.





\subsection{Partially Ordered Quantifiers}
Partially ordered (first-order) quantifiers were introduced by Henkin \cite{Hen61}, who noted that the usual first-order logic could not adequately express a situation in which for all $x$ there exists a $y$, and for all $z$ there is a $w$ such that $R(x, y, z, w)$, but the $y$ does not depend on $z$ and the $w$ does not depend on $x$. Since then, a number of authors have investigated such quantifiers and even infinitary versions of them \cite{Wal70}, \cite{End72}, \cite{Swe73}, etc. A comprehensive treatment was given in A.~K.~Swett's doctoral dissertation \cite{Swe73}. Swett assumed the partial order was well-founded. He gave an example to show that the semantics could be ill-defined if one dropped that assumption. It was shown by Walkoe \cite{Wal70} and Enderton \cite[p.~166]{End72} that sentences with finite partially ordered quantifiers are equivalent to existential second-order sentences, and hence by Corollary \ref{cor:MAExistential}  we obtain:

\begin{thm}
  Assuming $MA_{\omega_1}$, thin equivalence relations in first-order logic augmented with finite partially ordered quantifiers have $\leq\aleph_1$ equivalence classes.
\end{thm}

\begin{crl}
  Assuming $MA_{\omega_1}$, a countable theory in first-order logic augmented with finite partially ordered quantifiers has either a perfect set of pairwise non-isomorphic countable models or else has $\leq \aleph_1$ of them.
\end{crl}

We have not investigated the infinite partially ordered cases, except for the open and closed game quantifiers.

\section{Equivalence relations given by countable intersections of projective sets}\label{sec:Woodin}
The set of reals coding the countable models of a second-order \emph{sentence} can be arbitrarily high in the projective hierarchy.  The set of reals coding the countable models of a second-order \emph{theory} can thus be the intersection of countably many projective sets, and these projective sets may have unbounded projective complexity.  In order to extend our work to include second-order theories we therefore need information about the set of equivalence classes of an equivalence relation that is given by a countable intersection of projective sets. In this section we obtain results along these lines.  For this section we temporarily set aside our intended applications in second-order logic and focus only on the results we need in Descriptive Set Theory.  As such, we use the terminology of Descriptive Set Theory throughout; the only appearance of second-order logic in this section is at the very end, in Theorem \ref{thm:secondOrderTheory}.  In particular, we emphasize that in this section ``definable" is meant in the sense of Descriptive Set Theory.

At this point we should warn the reader that for this section not only do we assume familiarity with a substantial amount of Set Theory as we have implicitly so far but towards the end of the section we also assume some knowledge of \emph{Inner Model Theory}. For the set-theoretic background, i.e., up to and including the proof of Theorem \ref{thm:FM-infinitelymanyWoodins} modulo Lemmas \ref{lem:genericabsoluteness-infinitelymanyWoodins} and \ref{lem:uB-infinitelymanyWoodins}, standard textbooks such as \cite{Ka08, Je03, Sch14, Kun11} suffice. For the rest of this section we suppose familiarity with the basics of Inner Model Theory as for example covered in the first five sections of \cite{St10}.

Foreman and Magidor \cite{FM95} used a supercompact cardinal to produce a model of $\neg\CH$ in which every equivalence relation in $L(\bR)$ on the power set of $\bR$ has ${\leq}\aleph_1$ or a perfect set of inequivalent elements. Their model is the standard model for $\PFA$ but they could have instead used the forcing we apply in the proof of Theorem \ref{thm:FM-infinitelymanyWoodins}. By results of Woodin, they could have reduced the large cardinal hypothesis in their result from a supercompact cardinal to a proper class of Woodin cardinals or even to infinitely many Woodin cardinals with a measurable cardinal above them all.
There are two main consequences of supercompactness that Foreman-Magidor used in their proof. First, they used (their Fact 1.1) that sets in $L(\bR)$ are weakly homogeneously Suslin at the supercompact cardinal (see for example \cite[Problem 13.3]{Sch14} or \cite[Section 32]{Ka08} for the definition of $\kappa$-weakly homogeneously Suslin). Second, they used (their Theorem 1.2) generic absoluteness of $L(\bR)$ for forcings below the supercompact cardinal. See also the explanatory remarks after Theorem 1.2. For both of these applications, supercompactness is not necessary. Suppose $\kappa$ is an inaccessible cardinal and above $\kappa$ there are infinitely many Woodin cardinals with a measurable cardinal above them all. Then the Foreman-Magidor result follows using the forcing we use in the proof of Theorem \ref{thm:FM-infinitelymanyWoodins} (or by the forcing they used) by appealing to the following results of Woodin. First, generic absoluteness of $L(\bR)$ for forcings of size $\leq \kappa$ holds by \cite[Theorem 2.31]{Wo10}, see Theorem 6.1 and Remark 6.6 in \cite{St09DMT} for a proof. Second, every set of reals in $L(\bR)$ is $\kappa$-weakly homogeneously Suslin by \cite[Theorem 2.13]{Wo10}.

To emphasize, we state:

\begin{theorem}\label{thm:WoodinsWithMeasurableAbove}
    Suppose that there is an infinite sequence of Woodin cardinals with a measurable cardinal above all of them. Then it is consistent that $\neg\CH$ and every equivalence relation on $\mathbb{R}$ that is in $L(\mathbb{R})$ has $\leq \aleph_1$ or a perfect set of inequivalent elements. In particular, this holds for $\sigma$-projective equivalence relations.
\end{theorem}

For our application to variants of Morley's Theorem, we are interested in equivalence relations that are much simpler than arbitrary equivalence relations in $L(\bR)$. If we focus on equivalence relations given by countable intersections of projective sets we can reduce the large cardinal hypothesis to infinitely many Woodin cardinals (and even below, see the discussion below and Theorem \ref{thm:FM-laddermouse}). This leads to a significantly lower large cardinal hypothesis as the mice we consider do not have inner models with infinitely many Woodin cardinals. Defining mice is beyond the scope of this paper but we refer the interested reader to \cite{Sch14} or \cite{St10}.

\begin{theorem}\label{thm:FM-infinitelymanyWoodins}
Suppose there are infinitely many Woodin cardinals. Then there is a model of $\neg\CH$ in which every  equivalence relation on the power set of $\bR$ that is obtained as a countable intersection of projective sets has ${\leq}\aleph_1$ or a perfect set of inequivalent elements.
\end{theorem}

We recall the definition of Woodin cardinals for the reader's convenience. In contrast to other large cardinal notions such as measurable cardinals or supercompact cardinals, Woodin cardinals are not critical points of strong elementary embeddings. They are limits of such critical points in the following sense.

\begin{definition}
  \begin{enumerate}
      \item Let $\kappa < \delta$ be ordinals and $A \subseteq V_\delta$. Then $\kappa$ is called \emph{$A$-reflecting in $\delta$} if and only if for all $\eta < \delta$ there is an elementary embedding $i \colon V \rightarrow M$ with critical point $\kappa$ such that $i(\kappa) > \eta$ and \[ i(A) \cap V_\eta = A \cap V_\eta. \]
      \item A cardinal $\delta$ is a \emph{Woodin cardinal} if and only if for all $A \subseteq \delta$ there is some $\kappa < \delta$ that is $A$-reflecting in $\delta$.
  \end{enumerate}
\end{definition}

\begin{remark*}
 Note that the hypothesis of Theorem \ref{thm:FM-infinitelymanyWoodins} is not optimal. We will state a sharper theorem below but decided to start with this version as the statement of Theorem \ref{thm:FM-infinitelymanyWoodins} does not require any knowledge of Inner Model Theory.
\end{remark*}

The proof of Theorem \ref{thm:FM-infinitelymanyWoodins} uses the two following inner model theoretic lemmas that we prove below. Again, the large cardinal hypothesis in the statement of the lemmas is not optimal and we will improve it before proving the lemmas. We need one more definition before we can state the lemmas. 

Recall that projective sets of reals can be defined in second-order arithmetic with parameters by $\Sigma_n$ and $\Pi_n$ formulas, see for example \cite[Section 12, The Definability Context]{Ka08}.
We extend this hierarchy to formulas in the language $\cL_{\omega_1,\omega}$ defining $\sigma$-projective sets.  Extending the notion of equivalence of first-order formulas, we say two formulas $\phi(\overline{x}), \psi(\overline{x})$ of $\cL_{\omega_1, \omega}$ are \emph{equivalent} if $\mathcal{M} \models \forall \overline{x} (\phi(\overline{x}) \leftrightarrow \psi(\overline{x}))$ for every structure $\mathcal{M}$ in the relevant signature.  As in the projective case, we say a formula $\varphi$ is $\Sigma_{\alpha+1}$ for some $\alpha < \omega_1$ if it equivalent to $(\exists x) \psi$ for some $\Pi_\alpha$ formula $\psi$. For limit ordinals $\lambda<\omega_1$ we say a formula $\varphi$ is $\Sigma_{\lambda}$ if and only if it is equivalent to $\bigvee_{k<\omega} \psi_k$ for $\Sigma_{\lambda_k}$-formulas $\psi_k$ with $\lambda_k<\lambda$ for all $k<\omega$. 
Moreover, a formula $\varphi$ is $\Pi_{\alpha}$ for some $\alpha < \omega_1$ if it is equivalent to $\neg \psi$ for some $\Sigma_\alpha$ formula $\psi$. Below, we are interested in $\Pi_\omega$-formulas.

\begin{remark*}
 Note that a set of reals $A$ is a countable intersection of projective sets if and only if there is a $\Pi_\omega$-formula $\varphi$ and a parameter $z \in \BS$ such that $A$ is $\Pi_\omega(z)$ definable in second-order arithmetic in $z$, i.e., \[ A = \{ x \in \BS : \cA^2(z) \vDash \varphi(z) \}, \] where $\cA^2(z)$ is the two-sorted structure \[ (\omega, \BS, ap, +, \times, exp, <, 0,1,z) \] as in \cite[Section 12, The Definability Context]{Ka08}. Here $ap \colon {}^\omega\omega \rightarrow \omega$ denotes the binary operation of \emph{application}, i.e., $ap(x,n) = x(n)$.
\end{remark*}

\begin{lemma}\label{lem:genericabsoluteness-infinitelymanyWoodins}
Let $M$ be a model of $\ZFC$ with a sequence of Woodin cardinals $(\delta_i : i<\omega)$. Let $\bP$ be a partial order in $M$ with $|\bP| < \delta_0$ and let $g$ be $\bP$-generic over $M$. Then for every $\Pi_\omega$-formula $\varphi(v)$ and every $x \in (\BS)^M$, \[ (\cA^2(x))^M \vDash \varphi(x) \text{ if and only if } (\cA^2(x))^{M[g]} \vDash \varphi(x). \]
\end{lemma}

It is known that assuming large cardinals weaker than the existence of infinitely many Woodin cardinals (more precisely only the existence of certain inner models with finitely many Woodin cardinals) countable intersections of projective sets are determined, see for example \cite{AMS21}. It is not hard to obtain, from large cardinals, that they are ${<}\eta$-universally Baire for some ordinal $\eta$, as for example defined in \cite[Definition 32.21]{Je03} or \cite[Definition 8.6]{Sch14}. We recall the definition here. 

  \begin{definition}
    Let $(S,T)$ be trees on $\omega \times \kappa$ for some ordinal $\kappa$ and let $\eta$ be an ordinal. We say \emph{$(S,T)$ is $\eta$-absolutely complementing} if and only if \[ p[S] = \BS \setminus p[T] \] in every $\Col(\omega,\eta)$-generic extension of $V$.
  \end{definition}
  
  \begin{definition}[Feng-Magidor-Woodin, \cite{FMW92}]\label{def:uB} Let $A$ be a set of reals.
  \begin{enumerate}
      \item We say $A$ is \emph{${<}\eta$-universally Baire (${<}\eta$-uB)} if for every ordinal $\nu<\eta$, there are $\nu$-absolutely complementing trees $(S,T)$ with $p[S] = A$.
      \item We say $A$ is \emph{universally Baire (uB)} if it is ${<}\eta$-universally Baire for every ordinal $\eta$.
  \end{enumerate}
  \end{definition}

\begin{lemma}\label{lem:uB-infinitelymanyWoodins}
Let $M$ be a model of $\ZFC$ with a sequence of Woodin cardinals $(\delta_i : i<\omega)$. Then all countable intersections of projective sets in $M$ are ${<}\delta_0$-universally Baire.
\end{lemma}

We now prove Theorem \ref{thm:FM-infinitelymanyWoodins} modulo Lemmas \ref{lem:genericabsoluteness-infinitelymanyWoodins} and \ref{lem:uB-infinitelymanyWoodins} before we state the versions of these lemmas from weaker hypotheses.

\begin{proof}[Proof of Theorem \ref{thm:FM-infinitelymanyWoodins}]
Let $M$ be a model of $\ZFC$ with a sequence of Woodin cardinals $(\delta_i : i<\omega)$ and let $\kappa$ be the least inaccessible cardinal in $M$. In particular, $\kappa < \delta_0$. We consider a generic extension $M[G]$ of $M$ via the following forcing: Consider the countable support iteration $\bP_\kappa$ of length $\kappa$ of the following partial orders \[ \{ \dot{\bQ}_\alpha : \alpha < \kappa \}. \] At an even stage $\alpha<\kappa$ that is not a limit stage, let $\dot{\bQ}_\alpha$ be the usual countably closed collapse of the continuum (of the current stage of the iteration) to $\omega_1$. At limit stages $\alpha<\kappa$, let $\dot{\bQ}_\alpha$ be the trivial forcing.
At an odd stage $\beta<\kappa$ let $\dot{\bQ}_\beta$ add $\beta$ Cohen reals. Note that each individual forcing $\dot{\cQ}_\alpha$ and hence the whole iteration $\bP_\kappa$ is proper and in particular preserves $\aleph_1$. The final model $M[G]$ satisfies \[ 2^{\aleph_0} = \aleph_2 \] but $\CH$ holds at cofinally many initial segments $M[G \upharpoonright \alpha]$ of the iteration. Here $G \upharpoonright \alpha$ denotes the canonical restriction of the generic $G$ to the initial segment $\bQ_\alpha$ up to $\alpha$ of the iteration.

We claim that in $M[G]$ every equivalence relation on the power set of $\bR$ obtained by a countable intersection of projective sets has ${\leq}\aleph_1$ or a perfect set of inequivalent elements. Let $E$ be an equivalence relation on the power set of $\bR$ obtained by a countable intersection of projective sets and suppose $E$ is thin, i.e., it does not have a perfect set of inequivalent elements. We will argue that it has ${\leq}\aleph_1$ equivalence classes. 
Let $\varphi$ be a $\Pi_\omega$-formula defining $E$ in $M[G]$, i.e., \[ E = \{ (x,y) \in \bR^{M[G]} : (\cA^2(x,y))^{M[G]} \vDash \varphi(x,y) \}. \] Let $M[G \upharpoonright \alpha]$ be an initial segment of the iteration such that $\CH$ holds in $M[G \upharpoonright \alpha]$. Let $E_\alpha$ be the equivalence relation defined by \[ E_\alpha = \{ (x,y) \in \bR^{M[G \upharpoonright \alpha]} : (\cA^2(x,y))^{M[G \upharpoonright \alpha]} \vDash \varphi(x,y) \}. \]
Then, by $\CH$ in $M[G \upharpoonright \alpha]$, $E_\alpha$ has at most $\aleph_1$ equivalence classes in $M[G \upharpoonright \alpha]$. As the size of the forcing $\cP_\kappa$ is small, $M[G \upharpoonright \alpha]$ still has the sequence of Woodin cardinals $(\delta_i : i<\omega)$. So by Lemma \ref{lem:genericabsoluteness-infinitelymanyWoodins}, any pair of reals that is inequivalent under $E_\alpha$ is also $E$-inequivalent as $E$ is the equivalence relation given by $\varphi$.
The following claim finishes the proof of Theorem \ref{thm:FM-infinitelymanyWoodins} (modulo Lemmas \ref{lem:genericabsoluteness-infinitelymanyWoodins} and \ref{lem:uB-infinitelymanyWoodins}). 

\begin{claim}
Every real added by the forcing iteration from $M[G \upharpoonright \alpha]$ to $M[G]$ is $E$-equivalent to a real in $M[G \upharpoonright \alpha]$. In particular, $E$ also has at most $\aleph_1$ equivalence classes in $M[G]$.
\end{claim}
\textit{Proof.}
As in the previous claim, $M[G \upharpoonright \alpha]$ still has the sequence of Woodin cardinals $(\delta_i : i<\omega)$. So by Lemma \ref{lem:uB-infinitelymanyWoodins}, $E_\alpha$ is ${<}\delta_0$-universally Baire. This allows us to canonically extend $E_\alpha$ to equivalence relations in generic extensions for partial orders of size ${<}\delta_0$. Here the canonical extension of $E_\alpha$ to $M[G]$ is $E$.

Note that universal Baireness is all that is used of the hypothesis that the equivalence relation is weakly homogeneously Suslin in the proof of \cite[Theorem 3.4]{FM95}. Hence, we can apply this theorem to $E_\alpha$ and obtain that since $E_\alpha$ does not have a perfect set of inequivalent reals in $M[G]$, no real added by the forcing iteration from $M[G \upharpoonright \alpha]$ to $M[G]$ is $E$-equivalent to a real in $M[G \upharpoonright \alpha]$. 
\end{proof}

Now we turn to the proofs of Lemmas \ref{lem:genericabsoluteness-infinitelymanyWoodins} and \ref{lem:uB-infinitelymanyWoodins}. We in fact prove sharper versions of them that use a weaker hypothesis. Recall that the set of all countable models of a second-order $S$-theory for a countable signature $S$ is a countable intersection of projective sets. This is the case we are interested in here.

Before we can state the strengthenings of Lemmas \ref{lem:genericabsoluteness-infinitelymanyWoodins} and \ref{lem:uB-infinitelymanyWoodins} we need the following definition that introduces the large cardinal hypothesis we want to work from. Note that here we pass from large cardinals to inner models with large cardinals  --  a central theme in Inner Model Theory.

\begin{definition}
  We say $M$ is a \emph{ladder premouse above $\eta$} if and only if $M$ is a proper class fine structural premouse, as for example in \cite{St10}, with cardinals $\eta < \delta_0 < \delta_1 < \dots$ such that
  \begin{enumerate}
      \item $\eta$ is the second inaccessible cardinal in $M$,
      \item for each $n<\omega$, \[ M_n^\#(M|\delta_n) \models \text{``}\delta_n \text{ is Woodin''}, \] and
      \item for each $n<\omega$, \[ M_n^\#(M|\delta_n) \unlhd M. \]
  \end{enumerate}
\end{definition}

Here $M_n^\#(M|\delta_n)$ denotes the least $\omega_1$-iterable premouse constructed above $M|\delta_n$ that is not $n$-small above $\delta_n$, see for example \cite[p. 1660]{St10}, \cite{St95}, or the introduction of \cite{MSW} for the definition and some properties of this model.\footnote{Note that \cite[p. 1660]{St10} only mentions $\omega_1+1$-iterable models of the form $M_n^\#(X)$ for some $X$ but \cite{St95} and the introduction of \cite{MSW} deal more carefully with weaker concepts of iterability.} It is beyond the scope of this paper to formally introduce $M_n^\#(M|\delta_n)$ but we would like to mention that the premouse $M_n^\#(M|\delta_n)$ is a model that has $n$ Woodin cardinals above $\delta_n$. Nevertheless these Woodin cardinals will not remain Woodin in the full model $M$. So the minimal ladder premouse (if it exists) does not have Woodin cardinals. It has inner models with $n$ Woodin cardinals for every natural number $n$ but it does not have an inner model with infinitely many Woodin cardinals. This in particular explains why the \emph{consistency strength} of the existence of a ladder mouse is below the existence of a model with infinitely many Woodin cardinals as assumed in Theorem \ref{thm:FM-infinitelymanyWoodins}, even though there is no canonical large cardinal hypothesis below the existence of infinitely many Woodin cardinals that implies the existence of a ladder mouse. Note that the existence of a ladder premouse follows from the hypothesis of Theorem \ref{thm:FM-infinitelymanyWoodins}, i.e., in a model $M$ of $\ZFC$ with a sequence of Woodin cardinals $(\delta_i : i<\omega)$, by the techniques in \cite{St93}.

\begin{remark*}
A similar notion of ladder premouse figures prominently in the computation of the sets of reals of  initial segments of $L(\bR)$ as the sets of reals of mice, see for example \cite{Ru97} or the discussion in \cite[Section 8.4]{Uh16}.
\end{remark*}

Now we can state Lemmas \ref{lem:genericabsoluteness-infinitelymanyWoodins} and \ref{lem:uB-infinitelymanyWoodins} in terms of ladder premice. Note that if $M$ is a ladder mouse above $\eta$, then the reals $(\BS)^M$ of $M$ are clearly all contained in $M|\eta$, the model up to level $\eta$. 

\begin{lemma}\label{lem:genericabsoluteness-laddermouse}
Let $M$ be a ladder premouse above $\eta$ for some ordinal $\eta$. 
Let $\bP$ be a partial order in $M$ with $|\bP| < \eta$ and let $g$ be $\bP$-generic over $M$. Then for every $\Pi_\omega$ formula $\varphi(v)$ and every $x \in (\BS)^M$, \[ (\cA^2(x))^M \vDash \varphi(x) \text{ if and only if } (\cA^2(x))^{M[g]} \vDash \varphi(x). \]
\end{lemma}

\begin{lemma}\label{lem:uB-laddermouse}
Let $M$ be a ladder premouse above $\eta$ for some ordinal $\eta$ witnessed by cardinals $\delta_0 < \delta_1 < \dots$. 
Then all sets in $M$ that are obtained as a countable intersection of projective sets are ${<}\eta$-universally Baire.
\end{lemma}

For the proofs of Lemmas \ref{lem:genericabsoluteness-laddermouse} and \ref{lem:uB-laddermouse} we require that the reader is familiar with the basics of Inner Model Theory as for example covered in the first five sections of \cite{St10}.

\begin{proof}[Proof of Lemma \ref{lem:genericabsoluteness-laddermouse}]
Let $M$ be a ladder premouse above $\eta$ for some ordinal $\eta$. Then $M$ is clearly closed under the operations $z \mapsto M_n^\#(z)$ for all reals $z$ and natural numbers $n$. As the size of the forcing $\bP$ is below $\eta$, it is easy to see that any $\bP$-generic extension $M[G]$ of $M$ is also closed under the operations $z \mapsto M_n^\#(z)$ for all reals $z$ and natural numbers $n$. 
Therefore, the following claim finishes the proof of Lemma \ref{lem:genericabsoluteness-laddermouse}.

\begin{claim}
  Let $\varphi$ be a $\Sigma_\omega$ formula, say $\varphi = \bigvee_{k<\omega} \psi_k$, and let $N$ be any model that is closed under the operations $z \mapsto M_n^\#(z)$ for all reals $z$ and natural numbers $n$. Then for every $x \in (\BS)^N \cap (\BS)^V$, \[  (\cA^2(x))^N \vDash \varphi(x) \text{ if and only if } (\cA^2(x))^{V} \vDash \varphi(x). \]
\end{claim}
\textit{Proof.}
We start with the right-to-left implication. Let $x \in (\BS)^N \cap (\BS)^V$ and suppose \[ (\cA^2(x))^{V} \vDash \varphi(x). \] That means $(\cA^2(x))^{V} \vDash \psi_k(x)$ for some $k<\omega$. Suppose for notational simplicity that $\psi_k$ is a $\Sigma_k$ formula. Then, by correctness\footnote{This result is due to Woodin, see for example \cite[Lemma 1.17]{MSW} for the precise statement and a proof.} of the inner model $M_{k}^\#(x)$, \[ (\cA^2(x))^{M_{k}^\#(x)} \vDash \psi_k(x). \] By our assumption, $N$ is closed under the operation $z \mapsto M_k^\#(z)$ for all reals $z \in (\BS)^N$. In particular, $(M_k^\#(x))^N = (M_k^\#(x))^V$ and $(\cA^2(x))^{(M_{k}^\#(x))^N} \vDash \psi_k(x).$ Therefore, by correctness of $(M_k^\#(x))^N$ in $N$, \[ (\cA^2(x))^{N} \vDash \psi_k(x). \] Hence, $(\cA^2(x))^{N} \vDash \varphi(x)$, as desired. The left-to-right implication follows by the same argument.
\end{proof}

\begin{proof}[Proof of Lemma \ref{lem:uB-laddermouse}]
Let $A = \bigcap_{k<\omega} A_k$ be a countable intersection of projective sets $A_k$, $k<\omega$. Suppose for notational simplicity that each $A_k$ is a $\SIGMA^1_k$ set. Then by the usual argument for universal Baireness (or for obtaining homogeneously Suslin sets) from Woodin cardinals, see for example \cite[Problem 13.4]{Sch14}, \[ M_{k}^\#(M|\delta_{k}) \vDash A_k \text{ is } {<}\delta_k\text{-universally Baire}. \] In particular, $A_k$ is ${<}\eta$-universally Baire in $M_{k}^\#(M|\delta_k)$. Let $\nu < \eta$ and let $(S_k^\nu,T_k^\nu)$ be a pair of $\nu$-absolutely complementing trees with $p[S_k^\nu] = A_k$ in $M_{k}^\#(M|\delta_k)$ and in $M$ for each $k<\omega$. By combining the trees $(S_k^\nu,T_k^\nu)$ we can easily obtain a pair of $\nu$-absolutely complementing trees $(S^\nu,T^\nu)$ in $M$ with $p[S] = A$. Therefore, $A$ is ${<}\eta$-universally Baire in $M$.
\end{proof}

The proof of Theorem \ref{thm:FM-infinitelymanyWoodins} did not use any large cardinal hypothesis beyond an inaccessible cardinal except for the applications of Lemmas \ref{lem:genericabsoluteness-infinitelymanyWoodins} and \ref{lem:uB-infinitelymanyWoodins}. So by exchanging these with Lemmas \ref{lem:genericabsoluteness-laddermouse} and \ref{lem:uB-laddermouse} we obtain the following theorem.

\begin{theorem}\label{thm:FM-laddermouse}
Let $M$ be a ladder premouse above $\eta$ for some ordinal $\eta$. Then there is a model of $\neg\CH$ in which every equivalence relation on the power set of $\bR$ that is obtained as a countable intersection of projective sets has ${\leq}\aleph_1$ or a perfect set of inequivalent elements.
\end{theorem}

Finally, we apply the results of this section to second-order logic again.

\begin{thm}\label{thm:secondOrderTheory}
If there are infinitely many Woodin cardinals (or there exists a ladder premouse above some $\eta$), then there is a model of set theory in which the Absolute Morley Theorem holds for second-order theories in countable signatures.
\end{thm}
\begin{proof}
If $T$ is a second-order theory in a countable signature then the relation $\cong_T$ is obtained as an intersection of countably many projective sets, so this follows from Theorem \ref{thm:FM-infinitelymanyWoodins} (or Theorem \ref{thm:FM-laddermouse}).
\end{proof}


\bibliographystyle{amsplain}
\bibliography{b}

\providecommand{\bysame}{\leavevmode\hbox to3em{\hrulefill}\thinspace}
\providecommand{\MR}{\relax\ifhmode\unskip\space\fi MR }
\providecommand{\MRhref}[2]{%
  \href{http://www.ams.org/mathscinet-getitem?mr=#1}{#2}
}
\providecommand{\href}[2]{#2}
\begin{thebibliography}{10}

\bibitem{AMS21}
J.~P. Aguilera, S.~Müller, and P.~Schlicht, \emph{Long games and
  $\sigma$-projective sets}, Annals of Pure and Applied Logic \textbf{172}
  (2021), 102939.

\bibitem{BJ95}
T.~Bartoszynski and H.~Judah, \emph{Set theory, on the structure of the real
  line}, A. K. Peters Ltd., Wellesley, MA, 1995.

\bibitem{Todorcevic}
M.~Bekkali, \emph{Topics in set theory: Lebesgue measurability, large
  cardinals, forcing axioms, and $\rho$-functions; notes on lectures by {S}.
  {T}odorcevic}, Springer-Verlag, Berlin, 1991.

\bibitem{Burgess}
J.~Burgess, \emph{Infinitary languages and descriptive set theory}, Ph.D.
  thesis, University of California at Berkeley, Berkeley, California, 1974.

\bibitem{Burgess2}
J.~P. Burgess, \emph{Equivalences generated by families of {B}orel sets}, Proc.
  Amer. Math. Soc. \textbf{69} (1978), 323--326.

\bibitem{Davis}
M.~Davis, \emph{Infinite games of perfect information}, Advances in Game Theory
  (M.~Dresher, L.S. Shapley, and A.W. Tucker, eds.), Annals of Mathematics
  Studies, vol.~52, Princeton University Press, 1964, pp.~85--101.

\bibitem{Ebbinghaus}
E.-B. Ebbinghaus, \emph{Extended logics: {T}he general framework},
  Model-theoretic logics (J.~Barwise and S.~Feferman, eds.), Cambridge
  University Press, Cambridge, 2016, pp.~25--76.

\bibitem{End72}
H.~B. Enderton, \emph{A mathematical introduction to logic}, Academic Press,
  New York-London, 1972.

\bibitem{FMW92}
Q.~Feng, M.~Magidor, and H.~Woodin, \emph{{Universally Baire Sets of Reals}},
  Set Theory of the Continuum (New York, NY) (H.~Judah, W.~Just, and H.~Woodin,
  eds.), Springer US, 1992, pp.~203--242.

\bibitem{FM95}
M.~Foreman and M.~Magidor, \emph{{Large cardinals and definable counterexamples
  to the continuum hypothesis}}, Annals of Pure and Applied Logic \textbf{76}
  (1995), no.~1, 47--97.

\bibitem{Gao}
S.~Gao, \emph{Invariant descriptive set theory}, Chapman and Hall/CRC, New
  York, 2008.

\bibitem{HamkinsMO}
J.~D. Hamkins, \emph{At what level of the analytic hierarchy do {C}ohen reals
  lie?}, MathOverflow, URL: https://mathoverflow.net/q/186985 (version:
  2014-11-13)\\ Author URL:
  {https://mathoverflow.net/users/1946/joel-david-hamkins}.

\bibitem{Har77}
L.~Harrington, \emph{Long projective wellorderings}, Ann. Math. Logic
  \textbf{12} (1977), 1--24.

\bibitem{HSa}
L.~Harrington and R.~Sami, \emph{Equivalence relations, projective and beyond},
  Logic Colloquium’78, North-Holland, Amsterdam, 1979, pp.~247--264.

\bibitem{HSh}
L.~Harrington and S.~Shelah, \emph{Counting equivalence classes for
  co-$\kappa$-{S}ouslin equivalence relations}, Logic Colloquium’80, Prague,
  1980, North-Holland, Amsterdam, 1982, pp.~147--–152.

\bibitem{HartStarchenkoValeriote}
B.~Hart, S.~Starchenko, and M.~Valeriote, \emph{Vaught's conjecture for
  varieties}, Trans. Amer. Math. Soc. \textbf{342} (1994), no.~1, 173--196.

\bibitem{Hen61}
L.~Henkin, \emph{Some remarks on infinitely long formulas}, Infinitistic
  Methods (Proc. Sympos. Foundations of Math., Warsaw, 1959, Pergamon, Oxford;
  Pa\'nstwowe Wydawnictwo Naukowe, Warsaw, 1961, pp.~167--183.

\bibitem{Je03}
T.~J. Jech, \emph{Set theory}, Springer Monographs in Mathematics, Springer,
  2003.

\bibitem{Ka08}
A.~Kanamori, \emph{The higher infinite}, 2nd ed., Springer Monographs in
  Mathematics, Springer, 2009.

\bibitem{Kechris1978}
A.~Kechris, \emph{On transfinite sequences of projective sets with an
  application to {$\Sigma^1_2$}-equivalence relations}, Logic Colloquiuum '77
  (Amsterdam) (A.~Macintyre, L.~Pacholski, and J.~Paris, eds.), North-Holland,
  1978, pp.~155--160.

\bibitem{KechrisMoschovakis}
A.~S. Kechris, \emph{{AD} and projective ordinals}, Cabal seminar 76-77
  (Berlin) (A.S.~Kechris et. al., ed.), Lecture notes in mathematics, vol. 689,
  Springer-Verlag, 1978, pp.~91--132.

\bibitem{Kechris}
\bysame, \emph{Classical descriptive set theory}, Springer-Verlag, New York,
  1995.

\bibitem{Kun11}
K.~Kunen, \emph{Set theory}, College Publications, London, 2011.

\bibitem{Mar70}
D.~A. Martin, \emph{Measurable cardinals and analytic games}, Fund. Math.
  \textbf{66} (1970), 287--291.

\bibitem{Mar75}
\bysame, \emph{Borel determinacy}, Ann. Math. \textbf{102} (1973), 363--371.

\bibitem{Mar11}
\bysame, \emph{Projective sets and cardinal numbers: {S}ome questions related
  to the continuum problem}, The {C}abal Seminar, Volume {II}: {W}adge degrees
  and projective ordinals (A.S. Kechris, B.~L\"owe, and J.R. Steel, eds.),
  Lect. Notes Logic, no.~37, Cambridge University Press, 2011, pp.~484--508.

\bibitem{MS}
D.~A. Martin and R.~M. Solovay, \emph{Internal {C}ohen extensions}, Ann. Math.
  Logic \textbf{2} (1970), 143--178.

\bibitem{Martin1989}
D.~A. Martin and J.~Steel, \emph{A proof of projective determinacy}, J. Amer.
  Math. Soc. \textbf{2} (1989), no.~1, 71--125.

\bibitem{Mayer}
L.~Mayer, \emph{Vaught's conjecture for o-minimal theories}, J. Symb. Logic
  \textbf{53} (1988), no.~1, 146--159.

\bibitem{Morley1970}
M.~Morley, \emph{The number of countable models}, J. Symb. Logic \textbf{35}
  (1970), 14--18.

\bibitem{Mos72}
Y.~N. Moschovakis, \emph{The game quantifier}, Proc. Amer. Math. Soc.
  \textbf{31} (1972), 245--250.

\bibitem{Mos80}
\bysame, \emph{Descriptive set theory}, Mathematical Surveys and Monographs,
  vol. 155, American Mathematical Society, Providence, R.I., 2009.

\bibitem{MSW}
S.~M\"uller, R.~Schindler, and W.~H. Woodin, \emph{{Mice with finitely many
  Woodin cardinals from optimal determinacy hypotheses}}, Journal of
  Mathematical Logic \textbf{20} (2020).

\bibitem{Neeman2010}
I.~Neeman, \emph{Determinacy in {$L(\mathbb{R})$}}, Handbook of Set Theory, v.3
  (M.~Foreman and A.~Kanamori, eds.), Springer, Dordrecht, 2010,
  pp.~1877--1950.

\bibitem{Ru97}
M.~Rudominer, \emph{{Mouse sets}}, Annals of Pure and Applied Logic \textbf{87}
  (1997), no.~1, 1--100.

\bibitem{Sch14}
R.~Schindler, \emph{Set theory}, Universitext, Springer-Verlag, 2014.

\bibitem{Schlicht}
P.~Schlicht, \emph{Thin equivalence relations and inner models}, Ann. Pure and
  Appl. Logic \textbf{165} (2014), 1577--1625.

\bibitem{ShelahHarringtonMakkai}
S.~Shelah, L.~Harrington, and M.~Makkai, \emph{A proof of {V}aught's conjecture
  for $\omega$-stable theories}, Israel J. Math. \textbf{49} (1984), 259--280.

\bibitem{Sil80}
J.~H. Silver, \emph{Counting the number of equivalence classes of {B}orel and
  coanalytic equivalence relations}, Ann. Math. Logic \textbf{18} (1980),
  1--28.

\bibitem{Steel}
J.~Steel, \emph{On {V}aught's conjecture}, Cabal seminar 76-77 (Berlin)
  (A.S.~Kechris et. al., ed.), Lecture notes in mathematics, vol. 689,
  Springer-Verlag, 1978, pp.~193--208.

\bibitem{St93}
J.~R. Steel, \emph{Inner models with many {W}oodin cardinals}, Annals of Pure
  and Applied Logic \textbf{65} (1993), no.~2, 185--209.

\bibitem{St95}
\bysame, \emph{{Projectively well-ordered inner models}}, Annals of Pure and
  Applied Logic \textbf{74} (1995), 77--104.

\bibitem{St09DMT}
\bysame, \emph{{The derived model theorem}}, {Logic Colloquium 2006} (S.~B.
  Cooper, H.~Geuvers, A.~Pillay, and J.~Väänänen, eds.), {Lecture Notes in
  Logic}, Cambridge University Press, 2009, p.~280–327.

\bibitem{St10}
\bysame, \emph{{An Outline of Inner Model Theory}}, Handbook of Set Theory
  (M.~Foreman and A.~Kanamori, eds.), Springer, 2010.

\bibitem{Swe73}
A.~K. Swett, \emph{Herbrand's theorem in infinitary logic}, Ph.D. thesis,
  University of Toronto, 1973.

\bibitem{Uh16}
S.~{Uhlenbrock (now M\"uller)}, \emph{{Pure and Hybrid Mice with Finitely Many
  Woodin Cardinals from Levels of Determinacy}}, Ph.D. thesis, WWU Münster,
  2016.

\bibitem{Vaananen}
J.~V{\"a}{\"a}n{\"a}nen, \emph{Second-order and higher-order logic}, The
  Stanford Encyclopedia of Philosophy (E.~Zalta, ed.), The Metaphysics Research
  Lab Center for the Study of Language and Information at Stanford University,
  Stanford, 2020, URL:
  https://plato.stanford.edu/archives/fall2020/entries/logic-higher-order/.

\bibitem{Vaeaenaenen2023}
J.~V{\"a}{\"a}n{\"a}nen, \emph{Model theory of second order logic}, Beyond
  First Order Model Theory, II (J.~Iovino, ed.), CRC Press, Boca Raton, 2023.

\bibitem{Vaught}
R.~Vaught, \emph{Denumerable models of complete theories}, Proc. Sympos.
  Foundations of Mathematics, Infinitistic Methods (Warsaw), Pergamon Press,
  1961, pp.~303--321.

\bibitem{Velickovic}
B.~Velickovic, \emph{Forcing axioms and stationary sets}, Adv. Math.
  \textbf{94} (1992), 256--284.

\bibitem{Wal70}
W.~J. Walkoe, Jr., \emph{Finite partially-ordered quantification}, J. Symb.
  Logic \textbf{35} (1970), 535--555.

\bibitem{Wo10}
W.~H. Woodin, \emph{{The Axiom of Determinacy, Forcing Axioms, and the
  Nonstationary Ideal}}, {De Gruyter series in logic and its applications},
  vol.~1, De Gruyter, 2010.

\end{thebibliography}
\end{document}